\documentclass{amsart}
\usepackage{amssymb}
\newtheorem{Equation}{}[section]
\newtheorem{theorem}[Equation]{Theorem}
\newtheorem{proposition}[Equation]{Proposition}
\newtheorem{lemma}[Equation]{Lemma}
\newtheorem{corollary}[Equation]{Corollary}

\newtheorem{definition}[Equation]{Definition}

\newtheorem{remark}[Equation]{Remark}

\def\dd{\displaystyle}
\def\pa{\partial}
\def\cohom{\operatorname{H}}

\def\Hom{\operatorname{Hom}}
\def\End{\operatorname{End}}
\def\Ker{\operatorname{Ker}}
\def\Ind{\operatorname{Ind}}

\def\oH{\operatorname{H}}
\def\oK{\operatorname{K}}

\def\ch{\operatorname{ch}}

\def\Tr{\operatorname{Tr}}
\def\tr{\operatorname{tr}}

\def\C{\mathbb C}

\def\R{\mathbb R}
\def\Z{\mathbb Z}

\def\mfS{{\mathfrak S}}
\def\cA{{\mathcal A}}

\def\maC{{\mathcal C}}

\def\cG{{\mathcal G}}
\def\cH{{\mathcal H}}
\def\cO{{\mathcal O}}

\def\cS{{\mathcal S}}

\def\cU{{\mathcal U}}
\def\cV{{\mathcal V}}

\def\what{\widehat}
\def\wtit{\widetilde}
\def\tT{{\widetilde T}}

\def\te{{\tilde e}}

\def\wcA{\widetilde{\cA}}

\def\wL{\widetilde{L}}

\def\wT{\widetilde{T}}

\begin{document}



\title[Index theory and non-commutative geometry II \today]
{Index Theory and \\ Non-Commutative Geometry\\
II.  Dirac Operators and Index Bundles \\
\today}


\author[M-T. Benameur]{Moulay-Tahar Benameur}
\address{UMR 7122 du CNRS, Universit\'{e} de Metz, Ile du Saulcy, Metz, France}
\email{benameur@math.univ-metz.fr}
\author[J.  L.  Heitsch \today]{James L.  Heitsch}
\address{Mathematics, Statistics, and Computer Science, University of
Illinois at Chicago} \email{heitsch@math.uic.edu}

\begin{abstract}
When the index bundle of a longitudinal Dirac type operator is transversely smooth, we define its Chern character in Haefliger cohomology and relate it to the Chern character of the $K-$theory index. This result gives a concrete connection between the topology of the foliation and the longitudinal index formula. Moreover, the usual spectral assumption on the Novikov-Shubin invariants of the operator is improved. 
\end{abstract}
\maketitle
\tableofcontents

\section*{Introduction}
In this paper, we continue  our systematic study of the index theorem in Haefliger cohomology of foliations.  In    \cite{BH-I}, we defined  a Chern character for leafwise elliptic pseudodifferential operators on foliations.  By using Connes' extension in \cite{ConnesFundamental}, we then  translated the Connes-Skandalis $K-$theory index theorem \cite{CS:1984} into Haefliger cohomology, thus proving scalar index theorems in the presence of holonomy invariant currents. 

In order to get more insight into topological invariants of foliations, we extend here the results of \cite{Heitsch:1995} and
\cite{H-L:1999}, which tie the indices of a leafwise operator on a foliation of a compact manifold to the so-called index bundle of the operator.   In particular, we
show that for a  generalized Dirac operator $D$ along the leaves of a
foliation with Hausdorff graph, the Chern character of the analytic
index of $D$ coincides with the Chern character of the index
bundle of $D$.  As in \cite{Heitsch:1995} and \cite{H-L:1999}, we
assume that the projection onto the kernel of $D$ is transversely
smooth, and that the spectral projections of $D^{2}$ for the intervals
$(0,\epsilon)$ are transversely smooth, for $\epsilon$ sufficiently
small.  In those two papers, we assumed that the Novikov-Shubin
invariants of $D$ were greater than three times the codimension of
$F$.  Here we use the $\oK-$theory index and  we need only assume that they are greater than half the
codimension of $F$. More precisely, the pairings of these Chern characters with a given Haefliger  $2k-$current agree whenever the Novikov-Shubin invariants of $D$ are greater than $k$.  We conjecture that this theorem is still true
provided only that the Novikov-Shubin invariants are positive.  Note
that in the heat equation proof of the classical Atiyah-Singer
families index theorem, \cite{Bismut:1986}, it is assumed that there
is a uniform gap about zero in the spectrum of the operator, which
implies the conditions we assume on the spectral projections.

In \cite{Connes:1979, Connes:1981}, Connes extended the classical 
construction of Atiyah \cite{Atiyah:1975}  of the $L^2$ covering index theorem to leafwise elliptic operators on compact foliated
manifolds.  To do so he replaced the lifting and deck transformations used by Atiyah by a lifting
to the holonomy covers of the leaves invariant under the
natural action of the holonomy groupoid.  Moreover, he defined an
analytic index map from the K-theory of the tangent bundle of the
foliation to the K-theory of the C$^{*}$ algebra of the foliation,
which plays the role of the K-theory of the space of leaves.  In
\cite{CS:1984}, Connes and Skandalis defined a push forward map in
K-theory for any K-oriented map from a manifold to the space of leaves
of a foliation of a compact manifold.  This allowed them to define a
topological index map from the K-theory of the tangent bundle of the
foliation to the K-theory of its C$^{*}$ algebra.  Their main result
is that the analytic and topological index maps are equal, an
extension of the classical Atiyah-Singer families index theorem.
This theorem does not lead in general to a relation between the index of the operator and its  index bundle, by which we mean the (graded) projection  onto the kernel of the operator, even when this latter is transversely smooth so that its Chern character is well defined.  This index bundle, which lives in a von Neumann algebra of the foliation,  carries important information about the foliation, e.g. its Euler class, its higher signatures, etc. 

In this paper, we extend the Chern character to the index bundle of $D$, provided the projection onto the kernel of $D$ is transversely smooth.  Our main result is that, with the conditions given in the first paragraph, the Chern character of $D$ equals the Chern character of the index bundle of $D$. Since the Chern character of the index bundle equals the superconnection index defined in \cite{Heitsch:1995}, we obtain as a corollary the coincidence of the superconnection index with the Chern character of the analytic and topological indices.  This Chern character is readily computable and directly relates the index of $D$ with the topology of the foliation.  

 We point out the papers \cite{GL:2003, GL:2005} where Gorokhovsky and Lott prove, by a different method, an index theorem for longitudinal Dirac operators. 

Here is a brief outline of the paper.  In Section 1.,  we fix notation and briefly review some necessary material.  In Section 2., we extend our Chern character to the $K-$theory of the space of super-exponentially decaying operators on the leaves of a foliation, recall the construction of Dirac operators and the heat index idempotent. 
In Section 3., we review the construction of the Chern character we use, and extend it to the index bundle of a leafwise Dirac operator.  In Section 4., we prove our main theorem, Theorem \ref{limit}.  In Section 5., we show that the Chern character of the index bundle for $D$ defined here is the same as that defined in  \cite{Heitsch:1995} using Bismut superconnections.

Some of our results are valid for all foliations,  not just those with Hausfdorff groupoid.   We will alert the reader when we need to assume that the graph $\cG$ of $F$ is Hausdorff. It is also worth pointing out that  our results are valid if we replace the holonomy groupoid $\cG$ by any smooth groupoid between the monodromy and holonomy groupoids. 

{\em Achnowledgements.} The authors would like to thank A. Carey, A. Connes, M. Hilsum, E. Leichtnam, J. Lott, V. Nistor, and P. Piazza for many helpful discussions. We are especially indebted to G. Skandalis for many suggestions and remarks. 

\section{Notation and Review}\label{review}

Throughout this paper $M$ denotes a smooth compact Riemannian manifold of
dimension $n$, and $F$ denotes an oriented foliation of $M$ of
dimension $p$ and codimension $q$. So $n=p+q$.  The tangent bundle of $F$
will be denoted $TF$.
If $E \to N$ is a vector bundle over a manifold $N$, we denote the
space of smooth sections by $C^{\infty}(E)$ or by $C^{\infty}(N;E)$ if
we want to emphasize the base space of the bundle.  The compactly
supported sections are denoted by $C^{\infty}_{c}(E)$ or
$C^{\infty}_{c}(N;E)$.  The space of differential $k-$forms on $N$ is
denoted $\cA^{k}(N)$, and we set $\cA(N) = \oplus_{k\geq
0}\cA^{k}(N)$.  The space of compactly supported $k-$forms is denoted
$\cA_c^{k}(N)$, and $\cA_c(N) = \oplus_{k\geq 0}\cA_c^{k}(N)$.

The holonomy groupoid ${\mathcal G}$ of $F$ consists of equivalence
classes
of paths $\gamma:[0,1]\to M$ such that the image of $\gamma$ is
contained in a leaf of $F$.  Two such paths $\gamma_{1}$ and
$\gamma_{2}$  are
equivalent if $\gamma_{1}(0) = \gamma_{2}(0)$,
$\gamma_{1}(1) = \gamma_{2}(1)$, and the holonomy germ along
them is the same.  Two classes may be composed if the first ends where the
second begins, and the composition is just the juxtaposition of the
two
paths.  This makes $\cG$ a groupoid.  The space $\cG^{(0)}$ of units
of $\cG$ consists of the equivalence classes of the constant paths,
and we identify $\cG^{(0)}$ with $M$.

$\cG$ is a (in general non-Hausdorff) dimension $2p+q$ manifold.  The
basic
open sets defining its manifold structure are given as follows.  Let
${\cU}$ be a finite good cover of $M$ by foliation charts as defined in
\cite{H-L:1990}.  Given $U$ and $V$ in this cover and a leafwise path
$\gamma$ starting in $U$ and ending in $V$, define $(U, \gamma, V)$
to
be the set of equivalence classes of leafwise paths starting in $U$
and ending in $V$ which are homotopic to $\gamma$ through a homotopy
of leafwise paths whose end points remain in $U$ and $V$
respectively.
It is easy to see, using the holonomy defined by $\gamma$ from a
transversal in $U$ to a transversal in $V$, that if $U, V \simeq
\R^{p} \times \R^{q}$, then $(U,\gamma,V) \simeq \R^{p} \times
\R^{p} \times\R^{q}$.  If $\cG$ is non-Hausdorff, it is not true
that compact sets are always closed, nor that the closure of a compact
set is compact.  Because of this, we define the notion of having
compact support as follows.  Given a bundle $E$ over $\cG$ and any set
$(U, \gamma, V)$ as above, consider $E \,|\, (U, \gamma, V)$, the
restriction of $E$ to $(U, \gamma, V)$.  The space $C^{\infty}_{c}(E
\,|\, (U, \gamma, V))$ has a natural inclusion into the space of
sections of $E$ over $\cG$ by extending any element of
$C^{\infty}_{c}(E \,|\, (U, \gamma, V))$ to all of $\cG$ by defining
it to be zero outside $(U, \gamma, V)$.  We define the space
$C^{\infty}_{c}(E) = C^{\infty}_{c}(\cG;E) $ of smooth sections of $E$
over $\cG$ with compact support to be all finite sums $\sum s_{i}$
where each $s_i$ is the image of an element  in some
$C^{\infty}_{c}(E\,|\,(U, \gamma, V))$.  The space of smooth functions
with compact support on $\cG$, namely $C^{\infty}_{c}(\cG;\cG \times \R)$,
will
be denoted $C^{\infty}_{c}(\cG)$.  The metric on $M$ induces a
canonical metric on $\cG$, denoted $g_{0}$.  See \cite{Heitsch:1995}
for the construction.

The source and range maps of the groupoid $\cG$
are the two natural maps $s$, $r:\cG \to M$ given by
$s\bigl([\gamma]\bigr)=\gamma(0)$, $r\bigl([\gamma]\bigr)=\gamma(1)$.
$\cG$ has two natural transverse foliations $F_s$ and $F_r$ whose
leaves are respectively $\wL_x=s^{-1}(x)$, $\wL^x=r^{-1}(x)$ for $x\in
M$.  Note that $r:\wL_{x} \to L$ is the holonomy covering of $L$.

The Haefliger cohomology of $F$, \cite{Hae:1980}, is given as
follows. For each $U_i \in {\cU}$, let $T_i\subset U_i$ be a
transversal and set $T=\bigcup\,T_i$.  We may assume that the
closures
of the $T_i$ are disjoint.  Let $\cH$ be the holonomy pseudogroup
induced by $F$ on $T$.  Denote by $\cA^k_c(M/F)$ the quotient of
$\cA^k_c(T)$ by the vector subspace generated by elements of the form
$\alpha-h^*\alpha$ where $h\in \cH$ and $\alpha\in\cA^k_c(T)$ has
support contained in the range of $h$.  Give $\cA^k_c\bigl(M/F)$ the
quotient topology of the usual $C^\infty$ topology on $\cA^k_c(T)$,
so
this is not a Hausdorff space in general.  The exterior derivative
$d_T:\cA^k_c(T)\to \cA^{k+1}_c(T)$ induces a continuous differential
$d_H:\cA^k_c(M/F)\to \cA^{k+1}_c(M/F)$.  Note that $\cA^k_c(M/F)$ and
$d_H$ are independent of the choice of cover $\cU$.  The
complex $\{\cA^*_c(M/F),d_H\}$ and its cohomology $\oH^*_c(M/F)$ are,
respectively, the Haefliger forms and Haefliger cohomology of $F$.

As the bundle $TF$ is oriented, there is a continuous open surjective
linear map, called integration over the leaves,
$$
\int_F :\cA_c^{p+k}(M)\longrightarrow \cA^k_c(M/F)
$$
which commutes with the exterior derivatives $d_{M}$ and $d_{H}$.
Given $\omega \in \cA_c^{p+k}(M)$,  write $\omega =
\sum \omega_i$ where $\omega_i \in\cA^{p+k}_c(U_i)$.  Integrate $\omega_i$
along the fibers of the submersion $\pi_i:U_i\to T_i$ to obtain $\dd
\int_{U_i}
\, \omega_i  \in \cA^k_c(T_i)$.  Define $\dd \int_F \ \omega
\in\cA^k_c(M/F) $ to be the class of $\dd \sum_i \, \int_{U_i} \,
\omega_i$.  It is independent of the choice of the $\omega_i$
and of the cover $\cU$.  As $\dd \int_{F}$ commutes with $d_{M}$ and
$d_{H}$, it induces the map $\dd \int_F :\oH^{p+k}(M;\R) \to \oH^k_c(M/F)$.

\section{The $K-$theory index}\label{char}

In this section, we recall the definition of the analytic index of a  Dirac operator defined along the leaves of a foliation.   We begin with some general remarks about operators along the leaves of foliations.

Let $E_1$ and $E'_1$ be two complex vector bundles
over $M$ with Hermitian metrics and connections, and set $E = r^* E_1$ and $E' = r^* E'_1$
with the pulled back metrics and connections. 
A pseudo-differential $\cG$-operator with uniform support acting from $E$ to $E'$
is a {\em smooth} family $(P_x)_{x\in M}$ of $\cG$-invariant
pseudo-differential operators, where for each $x$, $P_x$ is an
operator  acting from $E \, | \,\wL_{x}$ to $E' \, | \,\wL_{x}$.  The
$\cG$-invariance property means
that for any $\gamma\in \wL_x^y=\wL_x \cap \wL^y$, we have $$
(\gamma\cdot P)_y = U_{\gamma} \circ P_x \circ U_{\gamma}^{-1} = P_y,
$$
where $U_{\gamma}$ denotes the operator on sections of any bundle induced by
the isomorphism $\gamma:\wL_{y} \to \wL_{x}$ given by
 composition with $\gamma$; for instance 
$$
U_{\gamma}: C_c^{\infty}(\wL_{x}, E) \longrightarrow
C_c^{\infty}(\wL_{y},E).
$$ 
The smoothness assumption is rigorously defined in
\cite{NistorWeinsteinXu}.  If we denote
by $K_x$ the Schwartz kernel of $P_x$, then the $\cG$-invariance
assumption implies that the family $(K_x)_{x\in M}$ induces a
distributional section $K$ of $\Hom(E, \widehat{E'})$ over $\cG$ which is smooth
outside  $\cG^{(0)} = M$.  Here $\widehat{E'}=s^*E'_1$, which is also the pullback bundle of $E'$ under the diffeomorphism $\gamma\mapsto \gamma^{-1}$. Since $M$ is compact, the
uniform support condition becomes the assumption that the support of
$K$ is compact in $\cG$.  The space of uniformly supported pseudo-differential $\cG-$operators
from $E$ to $E'$ is denoted $\Psi^{\infty}(\cG; E,E')$, and the space of uniformly supported
regularizing $\cG$-operators is denoted by $\Psi^{-\infty}(\cG;E,E')$. When $E'=E$ we simply denote the corresponding spaces by $\Psi^{\infty}(\cG; E)$ and $\Psi^{-\infty}(\cG;E)$.
The Schwartz Kernel Theorem identifies $\Psi^{-\infty}(\cG;E,E')$ with
$C_c^{\infty}(\cG,\Hom(E, \widehat{E'}))$, see \cite{Connes:1979,NistorWeinsteinXu}. 

 An element of
$\Psi^{\infty}(\cG; E,E')$ is elliptic if it is elliptic
when restricted to each leaf of $F_{s}$.  The parametrix theorem can be
extended to the foliated case and we have
\begin{proposition}\cite{Connes:1979}\ Let $P$ be a uniformly
supported elliptic pseudo-differential $\cG-$operator acting from $E$ to $E'$.  Then there
exists a uniformly supported pseudo-differential $\cG-$operator $Q$ acting from $E'$ to $E$
such that
$$
I_E - Q\circ P \in \Psi^{-\infty}(\cG; E) \text{ and } I_{E'} -
P\circ Q \in \Psi^{-\infty}(\cG; E').
$$
Here $I_E$ and $I_{E'}$ denote the identity operators of $E$ and $E'$ respectively.
\end{proposition}

A classical $K-$theory construction assigns to any uniformly supported elliptic pseudo-differential $\cG-$operator $P$ from $E$ to $E'$, a $K-$theory class
$$
\Ind_{a}(P) \in \oK_0(\Psi^{-\infty}(\cG; E\oplus E')) =
\oK_0(C_c^{\infty}(\cG,\Hom(E\oplus E')))
$$ 
called the analytic index of $P$ \cite{ConnesMoscovici,BH-I}. It will be useful to define this index class using functional calculus in a wider space of smoothing operators, so we now relax the uniform support condition and extend the above pseudodifferential calculus.

A super-exponentially decaying $\cG-$operator from  $E$ to $E'$ is a
family $P = (P_x)_{x\in M}$ of smoothing $\cG-$operators so
that its Schwartz kernel $P_{x}(y,z)$ is smooth in $x$, $y$, and $z$,
and satisfies
\begin{Equation}\label{Estsup}
Given non-negative integer multi indices $\alpha$, $\beta$,
and $\gamma$,  there are positive constants 
$\epsilon, C_{1}, \text{ and }C_{2}$, such that
for all $x \in M$, $y, z \in \wtit{L}_{x}$,
$$
\Vert \frac{\pa^{|\alpha |+|\beta | +|\gamma |} P_{x}(y,z)}{\pa x^\alpha \pa
y^\beta
\pa z^{\gamma}}\Vert \leq C_{1}\exp
\biggl[\frac{-d_x (y,z)^{1 + \epsilon}} {C_{2}}\biggr].$$
\end{Equation}
\noindent
Here $\pa/\pa x$ , $\pa/\pa y$, and $\pa/\pa z$ come from
coordinates obtained from a finite good cover $\cU$ of $M$ and $d_x(\, \,
,\, \,)$ is the distance on $\wtit
L_x$.  The space of all such operators is denoted
$\Psi^{-\infty}_{\mathfrak{S}}(\cG; E,E')$ or
$C_{\mathfrak{S}}^{\infty}({\cG}; \Hom(E,{\widehat{E'}}))$.  Again when $E'=E$ we denote the corresponding spaces by $\Psi^{-\infty}_{\mathfrak{S}}(\cG; E)$ and $C_{\mathfrak{S}}^{\infty}({\cG}; \Hom(E))$ for simplicity. When $E$ and $E'$ are trivial line bundles, we omit them and denote the corresponding spaces by $\Psi^{-\infty}_{\mathfrak{S}}(\cG)$ and $C_{\mathfrak{S}}^{\infty}({\cG})$.

\begin{lemma}
When $E'=E$, the space $\Psi^{-\infty}_{\mathfrak{S}}(\cG; E)$ is an algebra.
\end{lemma}

\begin{proof} Let $P$ and $Q \in \Psi^{-\infty}_{\mathfrak{S}}(\cG; E)$, with constants $\epsilon_{1}, C_{1}, C_{2}$ and $\epsilon_{2},D_{1}, D_{2}$
respectively, for the estimate given by Equation \ref{Estsup}.    We may replace 
$\epsilon_{1}$ and $\epsilon_{2}$ by $\epsilon = \min(\epsilon_{1},\epsilon_{2})$.  Set 
$\alpha = 1 + \epsilon$, $C= C_{1}D_{1}$ and $D = C_{2} + D_{2}$.  Then for $y, z \in \wL_{x}$,
$$
|P_{x}\circ Q_{x}(y,z)| = |\int_{\wL_{x}} P_{x}(y,w) Q_{x}(w,z) \, dw|
\leq \int_{\wL_{x}} C_{1} e^{-d(y,w)^{\alpha }/C_{2 }} D_{1} e^{-d(w,z)^{\alpha }/D_{2}}
\, dw
\leq
$$
$$
\int_{\wL_{x}} C e^{-d(y,w)^{\alpha }/D} e^{-d(w,z)^{\alpha }/D} \, dw =
C e^{-(d(y,z)^{\alpha }/2^{\alpha}D)}\int_{\wL_{x}} e^{-(d(y,w)^{\alpha } + d(w,z)^{\alpha } -
(d(y,z)/2)^{\alpha})/D} \, dw \leq
$$
$$
C e^{-(d(y,z)^{\alpha }/2^{\alpha}D)}\Big{[}\int_{S_{z}} e^{-d(y,w)^{\alpha }/D}  \, dw  +
\int_{S_{y}} e^{-d(w,z)^{\alpha }/D} \, dw \Big{]} \leq
$$
$$
C e^{-(d(y,z)^{\alpha }/2^{\alpha}D)}\Big{[}\int_{\wL_{x}} e^{-d(y,w)^{\alpha }/D}  \, dw  +
\int_{\wL_{x}} e^{-d(w,z)^{\alpha }/D} \, dw\Big{]},
$$
where
$$
S_{z} = \{ w \in \wL_{x} \, | \, d(w,z) \geq d(y,z)/2 \} \quad \text
{and} \quad S_{y} = \{ w \in \wL_{x} \, | \, d(y,w) \geq d(y,z)/2 \}.
$$
Now each of the integrals $\dd \int_{\wL_{x}} e^{-d(y,w)^{\alpha }/D}  \, dw$
and $\dd \int_{\wL_{x}} e^{-d(w,z)^{\alpha }/D} \, dw$ is bounded
independently of $x, y,$ and $z$.  This is a standard
argument for foliations of compact manifolds.  Since $M$ is compact,
the leaves $\wL_{x}$ have at most (uniformly bounded) exponential
growth, and the integrands are super-exponentially decaying with
uniform super-exponential bounds.  This gives us the estimate in
\ref{Estsup} for $P\circ Q$.

To get the estimate for the derivatives $\pa^{|\alpha |+|\beta | +|\gamma |}
(P\circ Q)_{x}(y,z)/\pa x^\alpha \pa y^\beta \pa z^{\gamma}$ we need
only note that these are finite sums of the form
$$
\sum_{\alpha_{1}+\alpha_{2} = \alpha} \int_{\wL_{x}}
\Big{(}\frac{\pa^{|\alpha_{1}|+|\beta|}P_{x}(y,w)}{ \pa x^{\alpha_{1}}
\pa y^{\beta}}\Big{)} \Big{(}\frac{\pa^{|\alpha_{2}|+|\gamma
|}Q_{x}(w,z)}{\pa x^{\alpha_{2}} \pa z^{\gamma}}\Big{)} \, dw.
$$
We can then repeat the argument above, using the estimates for the
individual integrands.

\end{proof}

There is a continuous embedding of algebras
$$
j_{\mathfrak{S}}: C_c^{\infty}({\cG}; \Hom(E\oplus E'))
\hookrightarrow C_{\mathfrak{S}}^{\infty}(\cG; \Hom(E\oplus E')),
$$
and we define the Schwartz analytic index $\Ind_a^{\mathfrak{S}}$ as
the composition of the analytic index $\Ind_a$ and the induced
morphism $j_{\mathfrak{S}*}:\oK_{*}({\cG}; \Hom(E\oplus E')) \to
\oK_0(C_{\mathfrak{S}}^{\infty}(\cG; \Hom( E\oplus E')))$.  So if
$P$ is a uniformly
supported elliptic pseudo-differential $\cG-$operator,
$$
\Ind_a^{\mathfrak{S}}(P) = j_{\mathfrak{S}*}(\Ind_a (P)) \ \ \in \ \
\oK_0(C_{\mathfrak{S}}^{\infty}(\cG; \Hom(E\oplus E'))).
$$
By classical arguments, see for instance \cite{MelroseNistor}, it is easy to check that $\Psi^{-\infty}_{\mathfrak{S}}(\cG; E,E')$ is a right module over the algebra $\Psi^{-\infty}_{\mathfrak{S}}(\cG)$.
The extended pseudodifferential calculus is defined by:
$$
\Psi^\infty_{\mathfrak{S}}(\cG; E,E') := \Psi^{-\infty}_{\mathfrak{S}}(\cG;  E,E') \otimes_{\Psi^{-\infty}(\cG)} \Psi^\infty (\cG;  E,E').
$$
It is generated by $\Psi^{\infty}(\cG; E,E')$ and
$\Psi^{-\infty}_{\mathfrak{S}}(\cG;  E,E')$.
When $E'=E$, we obtain in this way an algebra of pseudodifferential operators. The subspace
$\Psi^{-\infty}_{\mathfrak{S}}(\cG; E)$ is then an ideal
in the algebra $\Psi^{\infty}_{\mathfrak{S}}(\cG; E)$.  This is due to
the estimate given in \ref{Estsup}.  In particular, we may define
$\Ind_a^{\mathfrak{S}} (P)$ directly using a parametrix $Q \in
\Psi^{\infty}_{\mathfrak{S}}(\cG; E',E)$ and the classical
construction, and it is obvious that the two definitions agree. 

 The
construction of the Chern character $\ch_{a}:\oK_{0}(C_{c}^{\infty}(\cG;\Hom( E\oplus E')))\to \oH^{*}_{c}(M/F)$ in \cite{BH-I}, reviewed in Section \ref{chern} below, 
also extends to this case thanks
to Lemma 2.5, p.\ 443 of \cite{H-L:2002}.  Note that  this lemma requires one of the elements to be uniformly
exponentially decaying while the other must have uniformly bounded
coefficients.  But if an operator is uniformly exponentially decaying
it does have uniformly bounded coefficients.  Thus we have
$$
\ch_{a}^{\mathfrak{S}}:\oK_{0}(C_{\mathfrak{S}}^{\infty}(\cG;\Hom( E\oplus E')))  \longrightarrow
\oH^{*}_{c}(M/F)
$$
and
$$
\ch_{a}^{\mathfrak{S}} \circ j_{\mathfrak{S}*} = \ch_{a}.
$$
Finally, the formula for $\ch_{a}$ in Definition~\ref{ch} below also holds for
$\ch_{a}^{\mathfrak{S}}$.

Now assume that the dimension $p$ of $F$ is even and denote by $D$ a generalized Dirac operator for the foliation $F$.  One of the most important examples of such an operator is given by  the longitudinal Dirac operator with coefficients in a vector bundle over $M$.   It is defined as follows.  As above, let  $E_1$ be a complex vector bundle over $M$ with Hermitian metric and connection, and set $E = r^*(E_1)$ with the pulled back metric and connection. 
Assume that the tangent bundle $TF$ of $F$ is spin with a fixed spin structure.  Then $TF_{s}$
is also spin, and we endow it with the pulled back spin structure from $TF$. 
Denote by ${\cS} =
{\cS}^+\oplus {\cS}^{-}$ the bundle of spinors along the leaves of
$F_{s}$. Denote by $\nabla^0$ the connection on $TF_{s}$ given by the
orthogonal projection of the Levi-Civita connection for $g_0$ on
$T\cG$.  $\nabla^0$ in then the Levi-Civita connection on each leaf of
$F_{s}$ for the induced metric.  For all $x \in M$, $\nabla^0$ induces a connection
$\nabla^0$ on $\cS|{\wL_{x}}$ and we denote also by $\nabla^0$ the tensor
product connection on ${\cS} \otimes E |\wL_{x}$.  These data determine a
smooth family $D= \{D_x\}$ of  Dirac operators, where $D_x$
acts on sections of ${\cS \otimes E} |\wL_{x}$ as follows.  Let
$X_1,\ldots, X_p$ be a local oriented orthonormal basis of $T\wL_{x}$,
and set
$$
D_x=\sum^p_{i=1}\rho(X_i)\nabla^{0}_{X_i}
$$
where $\rho(X_i)$ is the Clifford action of $X_{i}$ on the bundle
${\cS} \otimes E |\wL_{x}$.  Then $D_x$ does not depend on the choice
of
the $X_{i}$, and it is an odd operator for the $\Z_{2}$ grading of
${\cS}\otimes E = ({\cS}^+ \otimes E)\oplus ({\cS}^{-}\otimes E)$.
Set
$D^{+} = D: C^{\infty}_{c}({\cS}^+ \otimes E) \to
C^{\infty}_{c}({\cS}^-
\otimes E)$
and $D^{-} = D: C^{\infty}_{c}({\cS}^- \otimes E) \to
C^{\infty}_{c}({\cS}^+
\otimes E)$. For more on generalized Dirac operators,   see \cite{LawsonMichelson}.

A super-exponentially decaying $\cG-$operator  on  $\cS \otimes E$ is defined to be an operator of the form
$$
A = \left(\begin{array}{cc} A_{11} & A_{12}\\ A_{21} &
A_{22}\end{array}\right),
$$
where each $A_{ij}$ is a smoothing operator  whose
Schwartz kernel $A_{ij,x}(y,z)$ is smooth in $x$, $y$, and $z$, and satisfies the estimate in
\ref{Estsup}.  $A_{11}$ maps sections of $\cS^{+} \otimes E$ to itself, $A_{12}$
maps sections of $\cS^{-}\otimes E$ to sections of $\cS^{+} \otimes
E$, etc.  The set  of all such operators is denoted $\Psi^{-\infty}_{\mathfrak{S}}(\cG; \cS \otimes E)$ or $C_{\mathfrak{S}}^{\infty}({\cG}; \Hom(\cS \otimes E))$.
If we unitalize  $\Psi^{-\infty}_{\mathfrak{S}}(\cG; \cS \otimes E)$ by adding two copies of $\C$ corresponding to the projections
$\pi_{\pm}:C^{\infty}_{c}(\cS \otimes E) \to C^{\infty}_{c}(\cS^{\pm}\otimes E)$, then we get a unital algebra that we denote by ${\widetilde\Psi}^{-\infty}_{\mathfrak{S}}(\cG; \cS \otimes E)$.  Note that
$\pi_{+} = \left(\begin{array}{cc} I & 0\\ 0 & 0 \end{array} \right) $ and
$\pi_{-} = \left(\begin{array}{cc} 0 & 0\\ 0 & I\end{array} \right) $.
Since the grading operator  $\alpha$ for $\cS = \cS^{+} \oplus \cS^{-}$ satisfies $\alpha= \pi_{+}  - \pi_{-}$,  $\alpha$ belongs to ${\widetilde\Psi}^{-\infty}_{\mathfrak{S}}(\cG; \cS \otimes E)$.

The odd operator $D$ is elliptic, so its analytic index is defined using a parametrix
$Q$ for $D$ which is also odd, i.e.
$$
Q =Q^{\pm}: C^{\infty}_{c}(\cS^{\pm} \otimes E) \longrightarrow
C^{\infty}_{c}(\cS^{\mp} \otimes E).
$$
Set
$$
S_{+} = I - Q^- \circ D^+ \text{ and } S_{-} = I - D^+ \circ Q^-
$$
so
$$
S_{\pm}: C^{\infty}_{c}(\cS^{\pm} \otimes E) \longrightarrow
C^{\infty}_{c}(\cS^{\pm} \otimes E).
$$
Using embeddings of our bundles in trivial bundles and computing the boundary map in $K-$theory, it is easy to see that the analytic index of $D$ is the $K-$theory class
\cite{ConnesMoscovici} in $\oK_0(\Psi^{-\infty}(\cG;\cS \otimes E)) =
\oK_{0}(C_c^{\infty}({\cG}; \Hom(\cS \otimes E))$,
$$
\Ind_a (D^+) = [{e}] - [\pi_-],
$$
where the idempotent $e$ is given by
\begin{Equation}\label{idem}  \hspace{1.5in}
$\dd {e} = \left(\begin{array}{cc} S_{+}^2 & - Q^- \circ (S_{-} + S_{-}^2)
\\
- S_{-} \circ D^+ & I - S_{-}^2 \end{array} \right).$
\end{Equation}

The class $[{e}] - [\pi_-]$ lives in the $K_0-$group of the unital algebra ${\widetilde \Psi}^{-\infty}(\cG;\cS \otimes E)$ but its image in the $K_0-$group of $\C\oplus \C$ under the map induced by
$$
p: \left(\begin{array}{cc} A_{11} + \lambda I_{\cS^{+} \otimes E} & A_{12}\\ A_{21} &
A_{22}+ \mu I_{\cS^{+} \otimes E}\end{array}\right) \longmapsto (\lambda, \mu),
$$
is trivial. Since this epimorphism admits a splitting homomorphism, it is clear that the kernel of the induced map $p_*$ is isomorphic  to the $K_0-$group of the non-unital algebra $\Psi^{-\infty}(\cG;\cS \otimes E)$. Hence, the index $\Ind_a (D^+) = [{e}] - [\pi_-]$ is well defined.

\begin{proposition}\label{Wassermann} \ 
Assume that $\cG$ is Hausdorff.  
Set
$$
P(tD) = \left[\begin{array}{cc} e^{-tD^-D^+} & (- e^{-tD^-D^+/2})
{\frac{I-e^{-tD^-D^+}}{t D^-D^+}} \sqrt{t}D^- \\
    &   \\
- e^{-tD^+D^-/2} \sqrt{t}D^+ & I - e^{-tD^+D^-} \end{array} \right].
$$
Then, for all $t >0$, $P(tD)$ is an idempotent in
${\widetilde\Psi}^{-\infty}_{\mathfrak{S}}(\cG; \cS \otimes E)$ and
$$
[P(tD)] - [\pi_-]= \Ind_a^{\mathfrak{S}}(D^+)\in
\oK_0({C}_{\mathfrak{S}}^{\infty}(\cG; \Hom(\cS \otimes E))).
$$
\end{proposition}

\begin{proof}\
It is classical that all the operators in $P(tD)$ (with the possible exception of the term  $\pi_-$) are smoothing when restricted to any $\wL_{x}$, so their Schwartz kernels are smooth when
restricted to any $\wL_{x}$.  Thus to check for smoothness, we need only check that
they are smooth transversely, i.e.  smooth in the variable $x \in M$.
The coefficients of the $D^{\pm}$ are smooth, and Corollary 3.11 of
\cite{Heitsch:1995} (which requires that $\cG$ be Hausdorff), says that the $e^{-tD^{\pm}D^{\mp}}$ are transversely smooth.
We will show presently that $\dd e^{-tD^-D^+/2} {\frac{I-e^{-tD^-D^+}}{t
D^-D^+}} \sqrt{t}D^- =   \sqrt{t}D^-e^{-tD^+D^-/2} {\frac{I-e^{-tD^+D^-}}{t D^+D^-}} $ is also transversely smooth.

By \cite{Heitsch:1995}, the Schwartz kernels $P_{t,x}^{\pm}(y,z)$ of the
$e^{-tD^{\pm} D^{\mp}}$ satisfy the following estimate.  Given a non-negative
integer $i$ and non-negative integer multi indices $\alpha$, $\beta$,
and $\gamma$, and a real number $T>0$, there is a constant $C >0$ such
that for all $x \in M$, $y, z \in \wtit{L}_{x}$, and $0 \leq t \leq
T$,
\begin{Equation}\label{super}
\hspace{1in} $\dd \Vert \frac{\pa^{i+|\alpha |+|\beta | +|\gamma |} P_{t,x}^{\pm}(y,z)}{\pa
t^{i}\pa x^\alpha \pa y^\beta \pa z^{\gamma}}\Vert \leq Ct^{-(p/2 +i
+|\alpha| + |\beta| + |\gamma|)} \exp\biggl[\frac{-d_x
(y,z)^2}{4t}\biggr].$
\end{Equation}
It follows immediately that the $e^{-tD^{\pm}D^{\mp}}$ and
the $e^{-tD^{\pm}D^{\mp}/2}$ satisfy the estimate in Equation \ref{Estsup},
and so also $e^{-tD^{+}D^{-}/2}\sqrt{t}D^{+}= \sqrt{t}D^{+}
e^{-tD^{-}D^{+}/2}$, since the derivatives of the
coefficients of $D^{\pm}$ are uniformly bounded on $\cG$.

To handle
$\dd e^{-tD^-D^+/2} {\frac{I-e^{-tD^-D^+}}{t D^-D^+}} \sqrt{t}D^- =
\sqrt{t}D^-e^{-tD^+D^-/2} {\frac{I-e^{-tD^+D^-}}{t D^+D^-}} $, note
that
$$
\frac{d}{ds} \Big{[}\frac{I-e^{-sD^+D^-}}{ D^+D^-}\Big{]} = e^{-sD^+D^-}, \quad \text{ so }
\quad \frac{I-e^{-tD^+D^-}}{t D^+D^-} = \frac{1}{t} \int_{0}^{t}
e^{-sD^+D^-} \,
ds.
$$
Thus
$$
\sqrt{t}D^-e^{-tD^+D^-/2} {\frac{I-e^{-tD^+D^-}}{t D^+D^-}}  =
\frac{\sqrt{t}D^{-}}{t}
\int_{0}^{t} e^{-(t/2+s)D^+D^-} \, ds = \frac{\sqrt{t}D^{-}}{t}
\int_{t/2}^{3t/2} e^{-sD^+D^-} \, ds.
$$
A simple calculation using Equation \ref{super} above then shows that for fixed $t$,
$\dd \sqrt{t}D^-e^{-tD^+D^-/2} {\frac{I-e^{-tD^+D^-}}{t D^+D^-}}$ is
transversely smooth and that it satisfies the estimate in Equation
\ref{Estsup}.

It is easy to check that the operator $Q(tD) =
Q^{\pm}(tD)$ where $Q^{-}(tD) = \dd {\frac{I-e^{-tD^-D^+/2}}{t D^-D^+}}
\sqrt{t}D^-$ and $Q^{+}(tD) = \dd {\frac{I-e^{-tD^+D^-/2}}{t
D^+D^-}}\sqrt{t}D^+$, is a parametrix for $\sqrt{t}D$.  The
corresponding idempotent $e$ given by Equation \ref{idem} is then $P(tD)$, so the
Schwartz analytic index of $tD$ is just $[P(tD)] - [\pi_-]$.   Since the index class only depends on the $K-$theory class of the principal symbol, it is clear that the $K-$theory class $[P(tD)] - [\pi_-]$ is independent of $t >0$.
\end{proof}

\section{The Chern character in Haefliger cohomology}\label{chern}

In this section we review the construction of the
Chern-Connes character in Haefliger cohomology  given in
\cite{BH-I}.  In view of our definition of the analytic index through the $K-$group of the unitalization ${\widetilde \Psi}^{-\infty}_{\mathfrak{S}}(\cG; \cS \otimes E)$, the Chern character is easy to express in terms of heat kernels. 
We may regard the connection $\nabla$ on $\cS\otimes E$ as an operator of degree one on $C^{\infty}(\cS \otimes E \otimes \land{T^{*}\cG})$ where on decomposable sections $\phi
\otimes \omega$, $\nabla(\phi \otimes \omega) = (\nabla\phi) \wedge
\omega + \phi \otimes d\omega.$ The foliation $F_s$ has normal bundle
$\nu_s^* \simeq s^{*}(T^{*}M)$, and $\nabla$ defines a {\em
quasi-connection} $\nabla^{\nu}$ acting on $C^{\infty}(\cS \otimes E \otimes
\land \nu_s^*)$ by the composition $$C^{\infty}(\cS \otimes E\otimes \land
\nu_s^*) \stackrel{i}{\longrightarrow} C^{\infty}(\cS \otimes  E\otimes
\land{T^{*}\cG})\stackrel{\nabla}{\longrightarrow} C^{\infty} (\cS \otimes E
\otimes \land{T^{*}\cG}) \stackrel{p_{\nu}}{\longrightarrow}
C^{\infty}(\cS \otimes E\otimes \land \nu_s^*),
$$
where $i$ is the inclusion and $p_{\nu}$ is induced by the
restriction $p_{\nu}:T^{*}{\cG} \to \nu^{*}_{s}$.

$C^{\infty}(\cS \otimes E\otimes \land \nu_s^*)$ is an $\cA(M)$-module where for
$\phi \in C^{\infty}(\cS \otimes E\otimes \land \nu_s^*)$, and $\omega \in
\cA(M)$, we set
$$
{\omega} \cdot \phi = p_{\nu}(s^{*}(\omega))\phi,
$$
where again the map $p_{\nu}:\cA({\cG}) \to
C^{\infty}(\land\nu_{s}^{*})$ is induced by the projection
$p_{\nu}:T^{*}{\cG} \to \nu^{*}_{s}$.

Recall $\Psi^{-\infty}(\cG;\cS \otimes E) \simeq C_c^{\infty}(\cG,\Hom(\cS \otimes E))$ the space of uniformly supported regularizing $\cG$-operators.  We may consider the algebra
$$
\Psi^{-\infty}({\cG}; \cS \otimes E) \widehat{\otimes} _{C^{\infty}(M)} \cA(M)
$$
as a subspace of the space of $\cA(M)$-equivariant endomorphisms of
$C^{\infty}(\cS \otimes E\otimes \land \nu_s^*)$ by using the $\cA(M)$ module
structure of $C^{\infty}(\cS \otimes E \otimes \wedge \nu^{*}_{s})$.

Denote by $\pa_{\nu}: \End(C^{\infty}(\cS \otimes E\otimes \land \nu_s^*)) \to
\End(C^{\infty}(\cS \otimes E\otimes \land \nu_s^*))$ the linear operator given by the graded commutator
$$
\pa_{\nu} (T) = [\nabla^{\nu}, T].
$$
The operator $\pa_{\nu}$ maps the space 
$$
\cA_c (\cG, \cS \otimes E):= \Psi^{-\infty}({\cG}; \cS \otimes E)
\widehat{\otimes}_{C^{\infty}(M)} \cA(M)
$$ 
to itself, and
$(\pa_{\nu})^2$ is given by the commutator with the curvature $\theta =
(\nabla^{\nu})^2$ of  $\nabla^{\nu}$.

In the same way, we consider the algebra 
$$
\cA_\mfS (\cG, \cS \otimes E):=\Psi^{-\infty}_\mfS ({\cG}; \cS \otimes E) \widehat{\otimes}_{C^{\infty}(M)} \cA(M)
$$
where $\Psi^{-\infty}_\mfS ({\cG}; \cS \otimes E)$ is the algebra of superexponentially decaying operators defined in the previous section. Then $\pa_\nu$ also acts on $\cA_\mfS (\cG, \cS \otimes E)$ with $(\pa_{\nu})^2$ given again by the commutator with the zero-th order differential operator $\theta$.

By the Schwartz kernel theorem, the algebra $\cA_\mfS (\cG, \cS \otimes E)$ is isomorphic to the algebra 
$$
C_\mfS^{\infty}({\cG};\Hom(\cS \otimes E))
\widehat{\otimes}_{C^{\infty}(M)} \cA(M).
$$
For any $T\in \cA_\mfS (\cG, \cS \otimes E)$,  define  the trace of $T$
to be the (compactly supported) Haefliger $k$-form $\Tr (T)$ given by
$$
\Tr(T) = \int_{F} \tr (K(\bar{x})) dx =  \int_{F} \tr (K(x,x)) dx,
$$
where $K$ is the smooth Schwartz kernel
of $T$, $\bar{x}$ is the class of the constant path at $x$,
$\tr(K(\bar{x}))$ is the usual trace of $K(\bar{x})
\in
\End((\cS \otimes E)_{\bar{x}}) \otimes \land T^* M_x$ and so belongs to
$\land T^* M_x$, and $dx$ is the leafwise volume form associated with
the fixed orientation of the foliation $F$.    The map 
$$
\Tr:\cA_\mfS (\cG,\cS \otimes  E) \longrightarrow \cA_{c}(M/F)
$$ 
is then a graded trace which satisfies 
$\Tr \circ \pa_{\nu} = d_{H} \circ \Tr$, see  \cite{BH-I}.

Since $\pa_{\nu} ^2$ is not necessarily zero, we used Connes' $X-$trick to construct a new graded differential algebra $({\wcA}_{\mfS }, \delta)$ out of the graded quasi-differential  algebra $(A_\mfS (\cG, \cS \otimes E), \pa_{\nu})$, see \cite{ConnesBook}, p.\ 229. 
First, note that  the curvature operator $\theta$ is a multiplier of $\cA_\mfS (\cG, \cS \otimes E)$.  As a 
vector space ${\wcA}_{\mfS } = M_2(A_\mfS (\cG,\cS \otimes  E))$.  An element 
${\widetilde T} =  \left(\begin{array}{cc}
{T}_{11} &{T}_{12}\\
{T}_{21} & {T}_{22} \end{array}\right)\in {\wcA}_{\mfS }$ is homogeneous of degree $\pa {\widetilde T} = k$ if
$$
k =  \pa { T}_{11} = \pa {T}_{12}+1 = \pa  {T}_{21} +1 = \pa  {T}_{22} +2.
$$
On homogeneous elements of ${\wcA}_{\mfS }$,  $\delta$ is given by
$$
\delta {\widetilde T} = \left(\begin{array}{cc}
\pa_{\nu} { T}_{11} & \pa_{\nu}{ T}_{12}\\
 - \pa_{\nu}{ T}_{21} & - \pa_{\nu}{ T}_{22}
\end{array}\right) 
+ 
\left(\begin{array}{cc} 0 & -\theta\\
1 & 0\end{array} \right){\widetilde T} 
+ 
(-1)^{\pa{\widetilde T}} {\widetilde T}
\left( \begin{array}{cc}
 0 & 1\\ -\theta & 0
\end{array}\right),
$$
and is extended to non-homogenous elements by linearity.
A  straightforward computation gives $\delta ^2=0$. 
For homogeneous  ${ T} \in A_\mfS (\cG,\cS \otimes  E)$, the differential $\delta$   on 
$ \left(\begin{array}{cc} {T}& 0\\
0 & 0\end{array}\right) \in {\wcA}_{\mfS }$ is given by
$$
\delta {\left(\begin{array}{cc} {T}& 0\\
0 & 0\end{array}\right)} = \left(\begin{array}{cc}
\pa_{\nu} { T} & (-1)^{\partial T}{T}\\
{ T} & 0 \end{array}\right).
$$

Set
$$
      \Theta=\left(
      \begin{array}{cc}
      1 &\quad 0\\
      0 & \quad \theta \end{array} \right)
$$
and define a new product on ${\wcA}_{\mfS }$ by
$$
      {\widetilde T} * {\widetilde T}' = {\widetilde T} \Theta {\widetilde T}'.
$$
This makes $({\wcA}_{\mfS }, \delta)$ a graded differential algebra.
For simplicity, we shall remove the multiplication $*$ from the notation and write $ {\widetilde T}  {\widetilde T}'$ for ${\widetilde T} * {\widetilde T}'$, when no confusion will occur.

The graded algebra $A_\mfS (\cG,\cS \otimes  E)$ embeds as a subalgebra of 
${\wcA}_{\mfS}$ by using the map
$$
	T \hookrightarrow \left( \begin{array}{cc} T & 0 \\ 0 & 0
	\end{array} \right).
$$
We shall therefore also denote by $T$ the image  in
${\wcA}_{\mfS}$ of any $T\in A_\mfS (\cG,\cS \otimes  E)$.

 For homogeneous ${\widetilde T}\in{\wcA}_{\mfS}$ define
$$
	\Phi ({\widetilde T}) = \Tr ({T}_{11}) - (-1)^{\pa {\widetilde
	T}} \Tr({T}_{22} \theta),
$$
and extend to arbitrary elements by linearity.
The map $\Phi:{\wcA}_{\mfS} \to \cA_c^*(M/F)$ is then a  graded trace,
and again we have $\Phi \circ \delta = d_H \circ \Phi$,  see  \cite{BH-I}.

The (algebraic) Chern-Connes character in the even case is the morphism 
$$
\ch_a:\oK_{0}(C_\mfS^{\infty}({\cG},\cS \otimes E)) = 
\oK_{0}(\Psi_\mfS^{-\infty}({\cG}; \cS \otimes E)) \longrightarrow \oH^*_{c}(M/F) 
$$
defined as follows.   
Denote  by ${\widehat{\Psi}}_\mfS^{-\infty}({\cG}; \cS \otimes E)$ the minimal unitalization of $\Psi_\mfS^{-\infty}({\cG}; \cS \otimes E)$.  This amounts to adding a copy of the complex numbers $\C$, so 
$$
{\widehat \Psi}^{-\infty}_{\mathfrak{S}}(\cG; \cS \otimes E) = \Psi^{-\infty}_{\mathfrak{S}}(\cG; \cS \otimes E) \oplus \C.
$$
Let
$M_{N}({\widehat{\Psi}}_\mfS^{-\infty}({\cG}; \cS \otimes E))$ be the space of $N \times N$
matrices with coefficients in  ${\widehat{\Psi}}_\mfS^{-\infty}({\cG}; \cS \otimes E)$.
Denote by 
$\tr:M_{N}(\Psi_\mfS^{-\infty}({\cG}; \cS\otimes E)) \to \Psi_\mfS^{-\infty}({\cG}; \cS \otimes E)$ the usual trace.

Recall the following:
\begin{theorem}\cite{BH-I}\label{BHI}
Let $B = [\te_{1}] - [\te_{2}]$ be an element of
$\oK_{0}(\Psi_\mfS^{-\infty}({\cG}; \cS \otimes E))$, where $\te_{1}=(e_1,\lambda_1)$
and $\te_{2}=(e_2,\lambda_2)$ are idempotents in
$M_{N}({\widehat{\Psi}}_\mfS^{-\infty}({\cG}; \cS \otimes E))$.  Then the Haefliger forms
$$
(\Phi\circ \tr)\Bigl{(}e_{1}\exp\left({\frac{-(\delta e_{1})^{2}}{2i \pi}}\right)\Bigr{)}\text{ and }
(\Phi\circ \tr) \Bigl{(}e_{2}\exp \left({\frac{-(\delta e_{2})^{2}}{2 i \pi}}\right)\Bigr{)}
$$
are closed and the Haefliger cohomology class of their difference depends only on $B$.

\end{theorem}

\begin{definition}\label{ch}   The algebraic Chern character $\ch_{a}(B)$
of $B$ is the Haefliger cohomology class
\begin{Equation}\label{ChernDef}
\hspace{1in} $ \dd  \ch_{a}(B) = \left[(\Phi\circ
\tr)\Bigl{(}e_{1}\exp\left({\frac{-(\delta
  e_{1})^{2}}{2i \pi}}\right)\Bigr{)}\right]- \left[(\Phi\circ \tr) \Bigl{(}e_{2}
\exp\left({\frac{-(\delta
  e_{2})^{2}}{2i\pi}}\right)\Bigr{)}\right].$
\end{Equation}
\end{definition}

In order to effectively compute the Chern character of the index of a generalized Dirac operator for $F$, we need some further results.
The exact sequence of algebras
$$
0\to \Psi^{-\infty}_{\mathfrak{S}}(\cG; \cS \otimes E) \stackrel{i}{\hookrightarrow} {\widetilde\Psi}^{-\infty}_{\mathfrak{S}}(\cG; \cS \otimes E) \stackrel{p}{\longrightarrow} \C^2 \to 0
$$
has a splitting homomorphism $\varrho: \C^2 \to {\widetilde\Psi}^{-\infty}_{\mathfrak{S}}(\cG; \cS \otimes E)$ given by $\varrho (\lambda,\mu) = \lambda\pi_+ + \mu \pi_-$. Therefore the kernel of the induced map
$$
p_*: K_0({\widetilde\Psi}^{-\infty}_{\mathfrak{S}}(\cG; \cS \otimes E)) \longrightarrow K_0(\C^2)\simeq \Z^2,
$$
is isomorphic to the group $K_0(\Psi^{-\infty}_{\mathfrak{S}}(\cG; \cS \otimes E))$. Denote by 
$p_0$ the obvious projection of ${\widehat \Psi}^{-\infty}_{\mathfrak{S}}(\cG; \cS \otimes E)$ onto $\C$. Then the inclusion map
$$
\beta : {\widehat \Psi}^{-\infty}_{\mathfrak{S}}(\cG; \cS \otimes E) \longrightarrow {\widetilde\Psi}^{-\infty}_{\mathfrak{S}}(\cG; \cS \otimes E),
$$
given by $\beta(T,\lambda)=T + \lambda \pi_+ +\lambda  \pi_-$ induces the isomorphism 
$$
\beta_*: K_0(\Psi^{-\infty}_{\mathfrak{S}}(\cG; \cS \otimes E) )=\Ker (p_{0,*}) \longrightarrow \Ker (p_{*}) \subset K_0({\widetilde\Psi}^{-\infty}_{\mathfrak{S}}(\cG; \cS \otimes E)).
$$

We shall use the universal graded algebra in the proof of Proposition \ref{idems} below, so we recall its definition. To any algebra $\maC$, there corresponds a (universal) differential graded algebra 
$\Omega(\maC)=\oplus_{n\geq 0} \Omega^n(\maC)$ which is defined by
$$
\Omega^0(\maC): = \maC \oplus \C, \text{ and for } n\geq 1, \Omega^n(\maC):= (\maC \oplus \C) \otimes \maC^{\otimes_n}.
$$
The differential $d:\Omega^n(\maC)\to \Omega^{n+1}(\maC)$ is defined for $a^j\in\maC$ and $c \in \C$ by
$$
	d\left[(a^0+c) \otimes a^1 \otimes \cdots \otimes a^n)\right]
	:= 1 \otimes a^0 \otimes a^1 \otimes \cdots \otimes a^n.
$$
It is clear that by definition $d^2=0$. The space $\Omega^n(\maC)$ is endowed with a natural right $\maC$-module structure (and hence right $\maC \oplus \C$-module structure) defined by
$$
	((a^0+c)  \otimes a^1\otimes \cdots \otimes a^n) a^{n+1} := (-1)^{n} \sum_{j=0}^n
	(-1)^j (a^0+c) \otimes \cdots \otimes a^ja^{j+1}\otimes \cdots
	\otimes a^{n+1}.
$$
The algebra structure of $\Omega(\maC)$ is defined by setting
$$
((a^0+c)  \otimes a^1\otimes \cdots \otimes a^n) (b^0 \otimes b^1 \otimes \cdots \otimes b^k) :=
 ((a^0+c)  \otimes a^1\otimes \cdots \otimes a^n)b^0 \otimes b^1 \otimes \cdots \otimes b^k 
$$
and
$$ 
((a^0+c)  \otimes a^1\otimes \cdots \otimes a^n)  c' \otimes b^1 \otimes \cdots \otimes b^k) :=
c' [(a^0+c) \otimes a^1\otimes \cdots \otimes a^n \otimes b^1 \otimes \cdots \otimes b^k].
 $$

A straightforward verification shows that $(\Omega(\maC) , d)$ is a differential graded algebra, see \cite{ConnesIHES}. We point out that by definition
$$
(a^0+c) da^1 \cdots da^n = (a^0+c) \otimes a^1\otimes \cdots \otimes a^n.
$$ 

The following  is known to experts. We give the proof for completeness, since it will be used in the sequel.

\begin{proposition}\label{idems}\
Let $\te$ and $\te'$ be two idempotents in $M_N({\widetilde\Psi}^{-\infty}_{\mathfrak{S}}(\cG; \cS \otimes E))$ such that $[\te] - [\te']$ belongs to the kernel of $p_*$. Then the Haefliger forms
\begin{multline*}
(\Phi\circ
\tr)\Bigl{(}(\te - (\varrho\circ p)(\te)) \exp\left({\frac{-(\delta
  (\te - (\varrho\circ p)(\te)))^{2}}{2i \pi}}\right)\Bigr{)} \text{ and } \\(\Phi\circ \tr) 
\Bigl{(}(\te' - (\varrho\circ p)(\te')) \exp\left({\frac{-(\delta  (\te' - (\varrho\circ p)(\te')))^{2}}{2i \pi}}\right)\Bigr{)}
\end{multline*}
are closed and we have the following equality in Haefliger cohomology:
\begin{multline*}
(\ch_a \circ \beta_*^{-1}) ([\te] - [\te']) = \left[(\Phi\circ
\tr)\Bigl{(}(\te - (\varrho\circ p)(\te)) \exp\left({\frac{-(\delta
  (\te - (\varrho\circ p)(\te)))^{2}}{2i \pi}}\right)\Bigr{)}\right]\\ - \left[(\Phi\circ \tr) 
\Bigl{(}(\te' - (\varrho\circ p)(\te')) \exp\left({\frac{-(\delta  (\te' - (\varrho\circ p)(\te')))^{2}}{2i \pi}}\right)\Bigr{)}\right]
\end{multline*}
\end{proposition}

\begin{proof}\
We define for every $k\geq 0$ a multilinear functional ${\widetilde \Phi}$ on the unital algebra
${\widetilde \Psi}^{-\infty}_\mfS ({\cG}; \cS \otimes E)$ by the equality
$$
{\widetilde \Phi} (\tT^0, \cdots, \tT^k) := \Phi( T^0 \delta T^1 \cdots \delta T^k) + \Phi (\delta (\Lambda^0 T^1)\delta T^2 \cdots \delta T^k),
$$
where  $\tT^j=T^j +\Lambda^j \in {\widetilde \Psi}^{-\infty}_\mfS ({\cG}; \cS \otimes E)$ with 
$$
T^j = \tT^j - (\varrho \circ p)(\tT^j) \in \Psi^{-\infty}_\mfS ({\cG}; \cS \otimes E)\text{ and } \Lambda^j=\varrho \circ p(\tT^j)=
\left(\begin{array}{cc} \lambda^j & 0 \\ 0 & \mu^j \end{array}\right) = \lambda^j \pi_+ + \mu^j \pi_-.
$$
Then ${\widetilde\Phi}$ is a functional on the universal differential graded algebra associated with ${\widetilde \Psi}^{-\infty}_\mfS ({\cG}; \cS \otimes E)$, see \cite{ConnesIHES} and also the bivariant constructions in \cite{CuntzQuillen,NistorInv}. More precisely, we set:
$$
{\widetilde \Phi} ((\tT^0+c) d\tT^1 \cdots  d\tT^k) := {\widetilde \Phi} (\tT^0, \cdots, \tT^k).
$$
We then have by definition
$$
({\widetilde \Phi} \circ d ) =0
$$
on the universal differential graded algebra associated with ${\widetilde \Psi}^{-\infty}_\mfS ({\cG}; \cS \otimes E)$.

For $\tT^j=T^j + \Lambda^j \in {\widetilde \Psi}^{-\infty}_\mfS ({\cG}; \cS \otimes E)$,
we have
\begin{eqnarray*}
(-1)^k{\widetilde \Phi}([\tT^0 d\tT^1 \cdots &&  \hspace{-0.33in}d\tT^k , \tT^{k+1}]) \quad  =  \quad(-1)^k{\widetilde \Phi}(\tT^0 d\tT^1 \cdots d\tT^k \tT^{k+1})
 - (-1)^k{\widetilde \Phi}( \tT^{k+1} \tT^0 d\tT^1 \cdots d\tT^k)\\ 
& \hspace{-0.3in} = & {\widetilde \Phi}(\tT^0\tT^1 d\tT^{2} \cdots d\tT^{k+1} ) 
+\sum_{j=1}^k (-1)^{j} {\widetilde \Phi}(\tT^0 d\tT^1 \cdots d\tT^{j-1} d(\tT^j\tT^{j+1}) d\tT^{j+2} \cdots d\tT^{k+1} ) \\ 
&-& (-1)^k{\widetilde \Phi} ( \tT^{k+1} \tT^0 d\tT^1 \cdots d\tT^k)\\ 
 & \hspace{-0.3in} =  & { \Phi}((T^0 T^1 +\Lambda^0 T^1 + T^0 \Lambda^1)\delta T^{2} \cdots \delta T^{k+1} ) +{\Phi}(\delta (\Lambda^0 \Lambda^1 T^2) \delta T^3 \cdots \delta T^{k+1} ) \\
& +&\sum_{j=1}^k (-1)^{j} \Phi(T^0 \delta T^1 \cdots \delta T^{j-1} \delta (T^jT^{j+1} + \Lambda^jT^{j+1} + T^j\Lambda^{j+1}) \delta T^{j+2} \cdots \delta T^{k+1} )\\ 
 & -& \Phi (\delta (\Lambda^0(T^1 T^2 + \Lambda^1 T^2 + T^1 \Lambda^2)) \delta T^3 \cdots \delta T^{k+1})\\
 & +& \sum_{j=2}^k (-1)^{j} \Phi (\delta (\Lambda^0T^1) \delta T^2 \cdots\delta T^{j-1}  \delta(T^jT^{j+1} + \Lambda^j T^{j+1} + T^j \Lambda^{j+1}) \delta T^{j+2} \cdots \delta T^{k+1})\\
 &-&  (-1)^k\Phi( ( T^{k+1}T^{0}+ T^{k+1}\Lambda^0+\Lambda^{k+1}T^0) \delta T^1 \cdots \delta T^k)\\
  &-& (-1)^k \Phi (\delta(\Lambda^{k+1} \Lambda^0 T^1) \delta T^2 \cdots \delta T^k).
\end{eqnarray*}

By using a connection which commutes with the grading we insure that $\pa^\nu (\Lambda) =0$ for any $\Lambda\in \C \pi_+ \oplus \C\pi_-$.   Thus, using the definitions of the product and the differential $\delta$, we can easily deduce the following relations for all $ \Lambda, T,  \Lambda'$, and $ T'$: 
\begin{Equation}\label{relations}
$\pa^\nu (\Lambda T) = \Lambda (\pa^\nu T),  \quad \pa^\nu (T \Lambda) = (\pa^\nu T)\Lambda,   \quad \theta \Lambda T = \Lambda \theta  T, \quad T \Lambda \delta (T') = T \delta (\Lambda T'), \quad \delta (T \Lambda) T' = (\delta T)  (\Lambda T'), $\\

\hspace{0.5in}$
\quad \delta(T\Lambda) \delta(T') = \delta (T) \delta (\Lambda T'),  \quad   T \Lambda \delta (\Lambda'T') = T \delta (\Lambda \Lambda'T') \quad \text{  and  } \quad 
\delta(TT') = \delta T T' + T \delta T' .
$
\end{Equation}
\noindent
It is then  a straightforward calculation that 
$$
{ \Phi}((T^0 T^1 +\Lambda^0 T^1 + T^0 \Lambda^1)\delta T^{2} \cdots \delta T^{k+1} )+
$$
$$
\sum_{j=1}^k (-1)^{j} \Phi(T^0 \delta T^1 \cdots \delta T^{j-1} \delta (T^jT^{j+1} + \Lambda^jT^{j+1} + T^j\Lambda^{j+1}) \delta T^{j+2} \cdots \delta T^{k+1} )
$$
collapses to 
$$
 \Phi (\Lambda^0T^1 \delta T^2 \cdots \delta T^{k+1}) +
 (-1)^k  \Phi(T^0 \delta T^1 \cdots \delta T^{k-1}( \delta T^k T^{k+1} + \delta(T^k\Lambda^{k+1}))), 
$$
and 
$$
{\Phi}(\delta (\Lambda^0 \Lambda^1 T^2) \delta T^3 \cdots \delta T^{k+1} )-
\Phi (\delta (\Lambda^0(T^1 T^2 + \Lambda^1 T^2 + T^1 \Lambda^2)) \delta T^3 \cdots \delta T^{k+1})+
$$
$$
\sum_{j=2}^k (-1)^{j} \Phi (\delta (\Lambda^0T^1) \delta T^2 \cdots\delta T^{j-1}  \delta(T^jT^{j+1} + \Lambda^j T^{j+1} + T^j \Lambda^{j+1}) \delta T^{j+2} \cdots \delta T^{k+1})
$$
collapses to 
$$
 -\Phi (\Lambda^0 T^1 \delta T^2  \cdots \delta T^{k+1}) +
 (-1)^k  \Phi(\delta (\Lambda^0 T^1) \delta T^2 \cdots \delta T^{k-1}( \delta T^k T^{k+1} + \delta(T^k \Lambda^{k+1}))).
$$
Substituting and multiplying by $(-1)^k$, we get
\begin{eqnarray*}
{\widetilde \Phi}([\tT^0 d\tT^1 \cdots d\tT^k , \tT^{k+1}]) &=&
(-1)^k \Phi (\Lambda^0T^1 \delta T^2 \cdots \delta T^{k+1}) \\
& + &  \Phi(T^0 \delta T^1 \cdots \delta T^k  T^{k+1}))\\
& +& \Phi(T^0 \delta T^1 \cdots \delta T^{k-1}\delta(T^k\Lambda^{k+1}))\\
& -&(-1)^k \Phi (\Lambda^0 T^1 \delta T^2  \cdots \delta T^{k+1}) \\
&+& \Phi(\delta (\Lambda^0 T^1) \delta T^2 \cdots \delta T^k T^{k+1})\\
& +&  \Phi(\delta (\Lambda^0 T^1) \delta T^2 \cdots \delta T^{k-1}\delta(T^k \Lambda^{k+1}))\\
&-& \Phi ( T^{k+1}T^{0}\delta T^1 \cdots \delta T^k)\\
&-& \Phi (  T^{k+1}\Lambda^0 \delta T^1 \cdots \delta T^k)\\
&-& \Phi ( \Lambda^{k+1}T^0 \delta T^1 \cdots \delta T^k)\\
&-&  \Phi (\delta(\Lambda^{k+1} \Lambda^0 T^1) \delta T^2 \cdots \delta T^k).
\end{eqnarray*}
The first and the fourth terms on the right cancel.  Using \ref{relations} and the trace property of $\Phi$ we have the following equations:
\begin{eqnarray*}
0 &=& \Phi (T^0 \delta T^1 \cdots \delta T^k T^{k+1}) 
-  \Phi(  T^{k+1}T^{0} \delta T^1 \cdots \delta T^k).\\
0 &=&  \Phi(T^0 \delta T^1 \cdots \delta T^{k-1}\delta(T^k\Lambda^{k+1}))
- \Phi(\Lambda^{k+1} T^0 \delta T^1 \cdots \delta T^k).\\
0 &=&  \Phi(\delta (\Lambda^0 T^1) \delta T^2 \cdots \delta T^k T^{k+1})
-  \Phi(  T^{k+1}\Lambda^0 \delta T^1 \cdots \delta T^k).\\
0 &=&  \Phi(\delta (\Lambda^0 T^1) \delta T^2 \cdots \delta T^{k-1}\delta(T^k \Lambda^{k+1}))
- \Phi (\delta(\Lambda^{k+1} \Lambda^0 T^1) \delta T^2 \cdots \delta T^k).
\end{eqnarray*}
Thus
$$
{\widetilde \Phi} ([\tT^0 d\tT^1 \cdots d\tT^k, \tT^{k+1}]) = 0.
$$

Now a classical argument shows that 
${\widetilde \Phi} $ is then a graded trace on the whole universal algebra associated with ${\widetilde \Psi}^{-\infty}_\mfS ({\cG}; \cS \otimes E)$.  

Given the above, we know that for any idempotent  $\te$ in the matrix algebra $M_N({\widetilde \Psi}^{-\infty}_\mfS ({\cG}; \cS \otimes E))$, the expression
$$
({\widetilde \Phi}\circ\tr)\Bigl{(}\te  \exp\left({\frac{-(d(\te)^2)}{2i \pi}}\right)\Bigr{)}
$$
 is a closed Haefliger form and that its cohomology class only depends on the $K-$theory class $[\te]$ of the idempotent $\te$, see for instance \cite{BH-I}. But note that this Haefliger differential form coincides up to exact Haefliger forms
with the differential form 
$$
(\Phi\circ
\tr)\Bigl{(}(\te - (\varrho\circ p)(\te)) \exp\left({\frac{-(\delta
  (\te - (\varrho\circ p)(\te)))^{2}}{2i \pi}}\right)\Bigr{)}
$$
which is then also closed and represents the same Haefliger cohomology class. Thus we deduce that the Haefliger class
\begin{multline*}
\left[(\Phi\circ
\tr)\Bigl{(}(\te - (\varrho\circ p)(\te)) \exp\left({\frac{-(\delta
  (\te - (\varrho\circ p)(\te)))^{2}}{2i \pi}}\right)\Bigr{)}\right] - \\ \left[(\Phi\circ \tr) 
\Bigl{(}(\te' - (\varrho\circ p)(\te')) \exp\left({\frac{-(\delta  (\te' - (\varrho\circ p)(\te')))^{2}}{2i \pi}}\right)\Bigr{)}\right],
\end{multline*}
is well defined and only depends on the $K-$theory class $[\te] - [\te']$. We denote it by ${\widetilde \ch}_a([\te] - [\te'])$. So we have the following morphism
$$
{\widetilde \ch}_a : K_0 ({\widetilde\Psi}^{-\infty}_{\mathfrak{S}}(\cG; \cS \otimes E)) \longrightarrow \cohom_c^*(M/F).
$$
The above construction applies also to the minimal unitalization ${\widehat\Psi}^{-\infty}_{\mathfrak{S}}(\cG; \cS \otimes E)$ of the algebra $\Psi^{-\infty}_{\mathfrak{S}}(\cG; \cS \otimes E)$ and yields a morphism 
$$
{\widehat \ch}_a:K_0({\widehat\Psi}^{-\infty}_{\mathfrak{S}}(\cG; \cS \otimes E)) \longrightarrow \cohom_c^*(M/F),
$$
whose restriction to  $K_0 (\Psi^{-\infty}_{\mathfrak{S}}(\cG; \cS \otimes E))$ is by definition the Chern character $\ch_a$. Note that ${\widehat \ch}_a$ is given by the same formula \eqref{ChernDef}, except that the $K-$theory element is no longer supposed to live in the kernel of 
$$
p_{0,*}:K_0({\widehat\Psi}^{-\infty}_{\mathfrak{S}}(\cG; \cS \otimes E)) \longrightarrow K_0(\Psi^{-\infty}_{\mathfrak{S}}(\cG; \cS \otimes E)).
$$
Now the map $\beta : {\widehat \Psi}^{-\infty}_{\mathfrak{S}}(\cG; \cS \otimes E)  \to  {\widetilde {\Psi}}^{-\infty}_{\mathfrak{S}}(\cG; \cS \otimes E)$  induces a well defined morphism of short exact sequences

\indent \hspace{0.2in}
\begin{picture}(415,80)
\put(0,60){$0$}
\put(0,10){$0$}
\put(10,64){\vector(1,0){20}}
\put(10,13){\vector(1,0){20}}
\put(40,60){$K_0(\Psi^{-\infty}_{\mathfrak{S}}(\cG; \cS \otimes E))$} \put(50,50){$\vector(0,-1){20}$}
\put(40,10){$K_0(\Psi^{-\infty}_{\mathfrak{S}}(\cG; \cS \otimes E))$} \put(65,35){$id$}

\put(150,68){$ i_{0,*}$}
\put(135,64){\vector(1,0){40}}
\put(135,13){\vector(1,0){40}}
\put(150,17){$i_{*}$}

\put(180,60){$K_0({\widehat \Psi^{-\infty}_{\mathfrak{S}}(\cG; \cS \otimes E)})$}
\put(210,50){ $\vector(0,-1){20}$}
\put(180,10){$K_0({\widetilde \Psi^{-\infty}_{\mathfrak{S}}(\cG; \cS \otimes E)})$}
\put(220,35){$\beta_*$}

\put(290,68){$p_{0,*}$}
\put(280,64){\vector(1,0){35}}
\put(280,13){\vector(1,0){35}}
\put(290,17){$p_{*}$}

\put(325,60){$K_0(\C)\simeq \Z$}
\put(345,50){ $\vector(0,-1){20}$}
\put(325,10){$K_0(\C^2)\simeq \Z^2$}
\put(355,35){$[\beta]_*$}

\put(385,64){\vector(1,0){20}}
\put(385,13){\vector(1,0){20}}

\put(415,60){$0$}
\put(415,10){$0.$}

\end{picture}

\noindent
Hence composing with ${\widetilde \ch}_a$ gives the following diagram which is commutative by the very definition of the maps: 

\indent \hspace{2in}

\begin{picture}(400,90)
\put(80,40){$0$}
\put(90,44){\vector(1,0){20}}
\put(115,40){$K_0(\Psi^{-\infty}_{\mathfrak{S}}(\cG; \cS \otimes E))$}

\put(175,70){$ i_{0,*}$}
\put(175,15){$i_{*}$}
\put(255,40){$\beta_*$}
\put(325,8){${\widetilde \ch}_a$}
\put(325,68){${\widehat \ch}_a$}

\put(180,55){\vector(1,1){25}}
\put(180,35){\vector(1,-1){25}}

\put(210,0){$K_0({\widetilde\Psi^{-\infty}_{\mathfrak{S}}(\cG; \cS \otimes E)})$} 
\put(210,80){$K_0({\widehat\Psi^{-\infty}_{\mathfrak{S}}(\cG; \cS \otimes E)})$}
\put(250,65){$\vector(0,-1){45}$} 
\put(310,5){$\vector(1,1){25}$} 
\put(310,75){$\vector(1,-1){25}$} 
\put(330,35 ){$\cohom_c(M/F)$.}
\end{picture}
\\
\\
In particular,  ${\widetilde \ch}_a \circ \beta_* = {\widehat \ch}_a$, so 
$$
{\widetilde \ch}_a \circ \beta_* \circ i_{0,*} = {\widehat \ch}_a \circ i_{0,*} =\ch_a.
$$
But, 
$$
\beta_* \circ i_{0,*}: K_0(\Psi^{-\infty}_{\mathfrak{S}}(\cG; \cS \otimes E))\longrightarrow \Ker p_*,
$$
is an  isomorphism, so we may define the Chern character directly on the group $K_0(\Psi^{-\infty}_{\mathfrak{S}}(\cG; \cS \otimes E)) =  \Ker p_*$. The proof is thus complete.
\end{proof}

\begin{corollary}\label{prop} Let $D$ be a  generalized Dirac operator for the foliation $F$ acting on the sections of the $\Z_2-$graded bundle $\cS\otimes E$. Let $P(tD)$ be the associated  idempotent in the algebra ${\widetilde\Psi}^{-\infty}_{\mathfrak{S}}(\cG; \cS \otimes E)$, as in Proposition \ref{Wassermann}.  Set 
$P_t  =  P(tD) -\pi_-$.  Then for all $t > 0$, the Haefliger form 
$$
 (\Phi\circ \tr)\Bigl{(}  P_t \exp \left[{\frac{-((\delta P_t)^2)}{2i\pi}}\right]  \Bigr{)},
$$
is closed and as Haefliger classes, we have the equality
$$
\ch_a (\Ind_a(D^+)) = \left[(\Phi\circ \tr)\Bigl{(}  P_t \exp \left[{\frac{-((\delta P_t)^2)}{2i\pi}}\right] \Bigr{)} \right].
$$
\end{corollary}

\begin{proof}\ The analytic $K-$theory index of $D$ in the $K-$theory group $K_0(\Psi^{-\infty}_{\mathfrak{S}}(\cG; \cS \otimes E))$ of superexeponentially decaying operators is given by
$$
\Ind_a(D^+) = [P(tD)] - [\pi_-] \quad \in \quad  \Ker \left( K_0({\widetilde\Psi}^{-\infty}_{\mathfrak{S}}(\cG; \cS \otimes E)) \to \Z^2 \right).
$$
Since the splitting map $\varrho: \C^2 \to {\widetilde\Psi}^{-\infty}_{\mathfrak{S}}(\cG; \cS \otimes E)$ is $\varrho (\lambda,\mu) = \lambda \pi_+ + \mu \pi_-,$
we have that 
$$
P(tD) - (\varrho\circ p) (P(tD)) = P_t \text{ and } \pi_- - (\varrho\circ p)(\pi_-)=0.
$$
Now apply Proposition \ref{idems}.
\end{proof}

In \cite{BH-I} we proved that the Chern character $\ch_a$ composed with the topological and analytic index maps of Connes-Skandalis \cite{CS:1984} yield the same map.  As a  particular case,  for any  generalized Dirac operator $D$ with coefficients in a Hermitian bundle $E_1$ over $M$, the Chern character of the topological index of $D$, denoted $\ch_{a}(\Ind_{t}(D^+))$, coincides with the 
Chern character of the analytic index of $D$, i.e.
$$
\ch_{a}(\Ind_{t}(D^+)) = \ch_{a}(\Ind_{a}(D^+)),
$$
and the common value of this Haefliger cohomology class is
$$
\ch_{a}(\Ind_{t}(D^+)) = \ch_{a}(\Ind_{a}(D^+)) = \int_F {\what A}(TF) \ch (E_{1}).
$$
Here ${\what A}(TF)$ is the usual ${\what A}$ genus of the tangent bundle of $F$, and $\ch$ is the usual Chern character of $E_{1}$.

In order to define the Chern character of the index bundle of $D$, we
need to assume that $P_{0}$, the projection onto the kernel of $D$, is
smooth.  Classical results imply that $P_{0}$ is smooth when
restricted to any leaf $\wL_{x}$, so what we are really assuming is that it is
transversely smooth.

Recall that  $\alpha =  \pi_{+} - \pi_{-}$ is the grading involution for ${\cS}\otimes E = ({\cS}^+\otimes E)\oplus ({\cS}^{-}\otimes E)$.  Then
$$
P_0 = \left[\begin{array}{cc} P^{+}_{0} & 0\\
0  & P^{-}_{0} \end{array} \right], \quad \text{ so   \quad   } \alpha P_0 = \left[\begin{array}{cc} P^{+}_{0} & 0\\
0  & - P^{-}_{0} \end{array} \right]
$$
is the super-projection onto the leafwise kernel of $D$, where
$P_{0}^{\pm}$ is projection onto the kernel of $D^{\pm}$.
Note  that $\pa_{\nu}\pi_{\pm} = 0$, provided we use a connection which preserves the splitting $\cS = \cS^{+} \oplus \cS^{-}$, which we assume that we do,
so  $\pa_{\nu} \alpha = 0$, and $\alpha \theta = \theta \alpha$.
Note also that $\alpha P_{0} = P_{0} \alpha$, so
$$
(\pa_{\nu}(\alpha P_{0}))^{2} =
\alpha^{2}(\pa_{\nu}P_{0})^{2} = (\pa_{\nu}P_{0})^{2} \text{  and  }
\alpha P_{0} \theta \alpha P_{0} = \alpha ^2  P_{0} \theta P_{0} =  P_{0} \theta P_{0}, 
\text{  which implies  } (\delta(\alpha P_0))^2 = (\delta P_0)^2.
$$

\begin{proposition}
The Haefliger form $\dd(\Phi\circ
\tr)\Bigl{(} \alpha P_{0}\exp({\frac{-((\delta
(\alpha P_{0}))^{2})}{2i\pi}})\Bigr{)} =    (\Phi\circ
\tr)\Bigl{(} \alpha P_{0}\exp({\frac{-((\delta P_{0})^{2})}{2i\pi}})\Bigr{)} $ is closed, and the Haefliger class it defines
depends only on $P_{0}$.
\end{proposition}

\begin{proof}\
Set $U=2P_0 -1$ then 
$$
\alpha U = U\alpha, U^2=I, UP_0 = P_0 = P_0 U \text{ and } U (\delta P_0) = \frac{1}{2} U (\delta U) = -\frac{1}{2} (\delta U) U = -(\delta P_0) U.
$$
Thus, for any $k\geq 0$, 
$$
(d_H \circ \Phi \circ \tr)\Bigl{(} \alpha P_{0} (\delta P_{0})^{2k})\Bigr{)} = 
(\Phi\circ\tr)\Bigl{(} \alpha (\delta P_{0})^{2k + 1} \Bigr{)}  = 
(\Phi\circ\tr)\Bigl{(}U^2 \alpha (\delta P_{0})^{2k + 1} \Bigr{)}.
$$
But,
$$
(\Phi\circ\tr)\Bigl{(}U^2 \alpha (\delta P_{0})^{2k + 1} \Bigr{)}=
( -1)^{2k+1} (\Phi\circ\tr)\Bigl{(}U \alpha  (\delta P_{0})^{2k + 1} U \Bigr{)} = 
- (\Phi\circ\tr) \Bigl{(}U^2 \alpha (\delta P_{0})^{2k + 1} \Bigr{)},
$$
so
$$
(d_H \circ \Phi \circ \tr)\Bigl{(} \alpha P_{0} (\delta P_{0})^{2k})\Bigr{)} = 0.
$$
In order to show the independence of the choice of connection, we use the relevant parts of the proof of Theorem 4.1 of \cite{BH-I}.  Indeed, it is obvious that the Poincar\'e argument developed there  still applies to the regularizing operator $P_0$ even though it may be non-compactly supported.   
\end{proof}

\begin{definition} The analytic Chern character $\ch_{a}(P_{0})$
of the index bundle of $D$ is the class of the Haefliger form
$\dd(\Phi\circ
\tr)\Bigl{(} \alpha P_{0}\exp({\frac{-((\delta
(\alpha P_{0}))^{2})}{2i\pi}})\Bigr{)} =    (\Phi\circ
\tr)\Bigl{(} \alpha P_{0}\exp({\frac{-((\delta P_{0})^{2})}{2i\pi}})\Bigr{)} $.
\end{definition}

Finally, an easy induction argument using the fact that for any idempotent $e$,
$e   (\pa_{\nu}e)^{2\ell - 1} e =0$ for all $\ell > 0$,
shows that 
$$
e (\delta e)^{2j} = 
 \left(\begin{array}{cc}
e \left(( \pa_{\nu}e)^{2}     +  e\theta e\right)^j   & 0\\
0 & 0
\end{array}\right). 
$$ 
Thus
\begin{Equation}\label{PO}
$\dd \hspace{1in} \ch_{a}(P_{0}) =)\Bigl{[} (\Tr\circ \tr)\Bigl{(}\alpha
P_{0}\exp(\frac{-((\pa_{\nu} P_{0})^{2} + P_{0}\theta P_{0})}{2i\pi})\Bigr{)}\Bigr{]}.$
\end{Equation}

\section{Proof of Main Theorem}\label{main}

Denote by $P_{\epsilon}$ the spectral projection for $D^{2}$ for the
interval $(0,\epsilon)$.  Recall that the Novikov-Shubin invariants
of
$D$ are greater than $k\geq 0$ provided that there is $\beta > k$ so
that
$$
(\Tr\circ \tr)(P_{\epsilon}) = (\Phi\circ \tr)(P_{\epsilon}) \mbox{  is  }  {\cO}(\epsilon^{\beta})   \mbox{  as  }   \epsilon \to 0.
$$
When we say a Haefliger form $\Psi$ depending on $\epsilon$ is
${{\cO}}(\epsilon^\beta)$ as $\epsilon \to 0$ we mean that there is a
constant $C>0$ so that the function on $T$, $\|\Psi\|_T \leq
C\epsilon^\beta$ as $\epsilon \to 0$.  Here $\| \quad \|_T$ is the
pointwise norm on forms on the transversal $T$ induced from the
metric on $M$.  To say that $\|\Psi\|_T \leq C\epsilon^\beta$, means that given any representative $\psi \in \Psi$, there is a constant $C_{\psi}$ so that
$\|\psi\|_T \leq C_{\psi} \epsilon^\beta$.   For a given cover, two representatives differ by a finite number of translations of local forms on transversals to other transversals so if this equation is satisfied for one representative with respect to a given cover, it is satisfied for all representatives with respect to that cover.  The fact that for any two good covers of $M$ by foliation
charts there is a integer $N$ so that any placque of the first cover
intersects at most $N$ placques of the second cover implies easily
that this condition does not depend on the choice of good cover.

We now prove our main theorem.

\begin{theorem}\label{limit}
Assume that $\cG$ is Hausdorff, and that the
Novikov-Shubin invariants of $D$ are greater than $q/2$.  Assume
further that the spectral projections $P_{0}$ and $P_{\epsilon}$ are
transversely smooth (for $\epsilon $ sufficiently small), and that $\pa_{\nu}P_{0}$ and $\pa_{\nu}P_{\epsilon}$ are bounded operators.
Then the analytic Chern character of the K-theory index of $D$ equals
the analytic Chern character of the index bundle of $D$, that is
$$
\ch_{a}(\Ind_{a}(D^+)) =  \ch_{a}([P_{0}]).
$$
\end{theorem}

Theorem \ref{limit} uses estimates on Novikov-Shubin invariants of $D$ to deduce the equality of the whole Chern character of the index bundle with that of the analytic index. We will actually prove  the following stronger theorem.

\begin{theorem}\label{limit-k}\
Assume again that $\cG$ is Hausdorff,  that the spectral projections $P_{0}$ and $P_{\epsilon}$ are transversely smooth (for $\epsilon $ sufficiently small), and that $\pa_{\nu}P_{0}$ and $\pa_{\nu}P_{\epsilon}$ are bounded operators.  For a fixed integer $k$ with $0\leq 2k\leq q$, assume that the Novikov-Shubin invariants of $D$ are greater than $k$.  
Then the $k^{th}$ component of the Chern character of the K-theory index of $D$ equals
the $k^{th}$ component of the Chern character of the index bundle of $D$, that is
$$
\ch^k_{a}(\Ind_{a}(D^+)) =  \ch^k_{a}([P_{0}]) \quad \in \quad \cohom_c^{2k}(M/F).
$$
\end{theorem}  

The proof of this theorem is rather long and involves a number of
complicated estimates.  For easier reading, we will split it into a
series of propositions and lemmas. Note that Theorem \ref{limit-k} implies Theorem \ref{limit}.

For the rest of this section, let $k$ be a fixed integer in the interval $[0,q/2]$.  By Corollary \ref{prop}, we need
only show that,
$$
\lim_{t \to \infty}(\Phi\circ \tr)\Bigl{(}  P_{t}(\delta P_{t})^{2k}\Bigr{)} = 
(\Phi\circ \tr)\Bigl{(}  \alpha P_{0}(\delta (\alpha P_{0}))^{2k}\Bigr{)} .
$$
If we ignore the minus signs in $P_{t}$, we see that the diagonal
terms give $e^{-t D^{2}}$, and the off diagonal terms are given by
$(P_{t})_{21} = (e^{-t D^{2}/2} \sqrt{t}D)_{21}$ and $\dd (P_{t})_{12} =
(e^{-t D^{2}/2}{\frac{I-e^{-t D^2}}{t D^2}}\sqrt{t} D )_{12}$.
Thus
$$
P_{t} = \pi_{+}e^{-tD^{2}}\pi_{+} - \pi_{-}e^{-tD^{2}}\pi_{-} -
\pi_{-}e^{-t D^{2}/2} \sqrt{t}D\pi_{+} - \pi_{+}e^{-t
D^{2}/2}{\frac{I-e^{-t D^2}}{t D^2}}\sqrt{t} D \pi_{-}.
$$
As the connection $\nabla$ used in the definition of $\pa_{\nu}$
preserves the splitting ${\cS}\otimes E = ({\cS}^+ \otimes E)\oplus
({\cS}^{-}\otimes E)$,  $\pa_{\nu} \pi_{\pm} = 0$, and we may work
with the operators $e^{-tD^{2}}$, $e^{-t D^{2}/2} \sqrt{t}D$, and
$\dd e^{-t D^{2}/2}{\frac{I-e^{-t D^2}}{t D^2}}\sqrt{t} D$ in what follows
instead of the (more notationally complicated) entries of $P_{t}$.

We will assume that the reader is familiar with the Spectral Mapping
Theorem, see for instance \cite{Reed-Simon}, and how to use it to
compute bounds on norms, strong convergence, etc.  This theorem gives
that
for $\ell \geq 0$, the norms of the operators $D^{\ell}e^{-t D^{2}}$,
$D^{\ell}e^{-t D^{2}/2}\sqrt{t}D$ and $\dd D^{\ell}e^{-t
D^{2}/2}{\frac{I-e^{-t D^2}}{t D^2}}\sqrt{t}D$ are uniformly bounded
as $t \to \infty$.  In addition, as $t \to \infty$, $D^{\ell}e^{-t
D^{2}/2}\sqrt{t}D$ and $\dd D^{\ell}e^{-t D^{2}/2}{\frac{I-e^{-t
D^2}}{tD^2}}\sqrt{t}D$ converge in norm to zero for $\ell \geq 0$,
and
for $\ell > 0$, $D^{\ell}e^{-t D^{2}}$ also converges in norm to
zero. 

Choose $\delta$ so that
$$
-1 < \delta < \frac{-k}{\beta} < 0
$$
and couple $\epsilon$ to $t$ by setting
$$
\epsilon = t^{\delta}.
$$

Because of the uniformly bounded geometry of the leaves of
$F_{s}$, which follows from the fact that all the structures we use on
$\cG$ are pulled back from the compact manifold $M$, the leafwise
estimates we give below are uniform over all leaves of $F_{s}$.

Denote by $Q_{\epsilon}$ the spectral projection for $D^{2}$ for the
interval
$[\epsilon,\infty)$.  Since $I = P_{0} + P_{\epsilon} +
Q_{\epsilon}$, the operator $\pa_{\nu}Q_{\epsilon}$ is bounded.
Now consider 
$$P_{t} = P_{0}P_{t}P_{0} +
P_{\epsilon}P_{t}P_{\epsilon} + Q_{\epsilon}P_{t}Q_{\epsilon} = \alpha P_0
+ P_{\epsilon}P_{t}P_{\epsilon} + Q_{\epsilon}P_{t}Q_{\epsilon}.
$$

\begin{proposition}\label{estimates} As $t \to \infty$,\\
\indent (i) $||Q_{\epsilon}P_{t}Q_{\epsilon}||$ is
bounded by a multiple of $e^{-(t^{(1+ \delta)}/32)}$,\\
\indent(ii) $||\pa_{\nu}(Q_{\epsilon}P_{t}Q_{\epsilon})||$ is
bounded by a multiple of $e^{-(t^{(1+ \delta)}/32)}$,\\
\indent (iii)  $||P_{\epsilon}P_{t}P_{\epsilon}||$ is
bounded,\\
\indent (iv) $||\pa_{\nu}(P_{\epsilon}P_{t}P_{\epsilon})||$ is
bounded by a multiple of $t^{(\frac{1}{2}+a)}$, for any $a>0$.

\end{proposition}

\begin{remark}\
The coefficient $\frac{1}{32}$ in (i) and (ii) can be improved very easily but this does not allow us to improve the assumption on the Novikov-Shubin invariants.
\end{remark}

\begin{proof}  
Note that the element
$$\pa_{\nu}( Q_{\epsilon}P_{t}Q_{\epsilon}) =
\pa_{\nu}(Q_{\epsilon})P_{t}Q_{\epsilon}+
Q_{\epsilon}P_{t}\pa_{\nu}(Q_{\epsilon})+
Q_{\epsilon}\pa_{\nu}(P_{t})Q_{\epsilon}
$$
and $||\pa_{\nu}(Q_{\epsilon})||$ is bounded.
We may write $P_{t} = e^{-tD^{2}/4}  \what{P}_{t} = \what{P}_{t} e^{-tD^{2}/4} $
where
$$
\what{P}_{t} = \left[\begin{array}{cc} e^{-3tD^-D^+/4} & (-
e^{-tD^-D^+/4}) {\dd \frac{I-e^{-tD^-D^+}}{t D^-D^+}} \sqrt{t}D^- \\
    &   \\
- e^{-tD^+D^-/4} \sqrt{t}D^+ & - e^{-3tD^+D^-/4} \end{array} \right].
$$
$\what{P}_{t}$ has essentially the same properties as $P_{t}$, in
particular its norm is bounded independently of $t$.  Since
$||e^{-tD^{2}/4}Q_{\epsilon}|| = ||Q_{\epsilon}e^{-tD^{2}/4}||\leq
e^{-t\epsilon/4} = e^{(-t^{(1+ \delta)}/4)}$, we have that
$||P_{t}Q_{\epsilon}||$ and $||Q_{\epsilon}P_{t}||$ (so also
$||Q_{\epsilon}P_{t}Q_{\epsilon}||$,
$||\pa_{\nu}(Q_{\epsilon})P_{t}Q_{\epsilon}||$ and
$||Q_{\epsilon}P_{t}\pa_{\nu}(Q_{\epsilon})||$) are bounded by a
multiple of $e^{(-t^{(1+ \delta)}/4)}$.    Thus we have $(i)$ of the Proposition, and to establish $(ii)$  we need only consider
the term $Q_{\epsilon}\pa_{\nu}(P_{t})Q_{\epsilon}$.

\begin{lemma}\label{estk}
$||Q_{\epsilon}\pa_{\nu}(e^{-tD^{2}/k})Q_{\epsilon}||$ is bounded by
a
multiple of  $e^{-(t^{1+\delta}/8k)}$.
\end{lemma}

\begin{proof}
Recall the foliation Duhamel formula of \cite{Heitsch:1995} (which requires that $\cG$ be Hausdorff) which
states that
$$
\pa_{\nu}(e^{-t D^{2}}) =
-\int^t_0 e^{-s D^{2}} \partial_{\nu}
(D^{2}) e^{(s-t) D^{2}} ds.
$$
Thus
$$
Q_{\epsilon}\pa_{\nu}(e^{-t D^{2}})Q_{\epsilon} =
-\int^t_0 Q_{\epsilon} e^{-s D^{2}} \partial_{\nu}
(D^{2}) e^{(s-t) D^{2}} Q_{\epsilon}ds =
$$
$$
-\int^t_{t/2} Q_{\epsilon} e^{-s D^{2}} \partial_{\nu}
(D^{2}) e^{(s-t) D^{2}} Q_{\epsilon}ds  -\int^{t/2}_0 Q_{\epsilon}
e^{-s
D^{2}} \partial_{\nu} (D^{2}) e^{(s-t) D^{2}} Q_{\epsilon}ds.
$$
The norm of the first integral satisfies

\noindent $\begin{array}{lll} \dd ||\int^t_{t/2} Q_{\epsilon}
e^{-sD^{2}} \partial_{\nu} (D^{2}) e^{(s-t) D^{2}} Q_{\epsilon}ds|| &
\leq
&\dd \int^t_{t/2} ||Q_{\epsilon} e^{-s D^{2}} \partial_{\nu} (D^{2})
e^{(s-t) D^{2}} Q_{\epsilon}||ds\\
&&\\
& \leq &
\dd \int^t_{t/2} ||Q_{\epsilon} e^{-\frac{s}{2} D^{2}}||\,
||e^{-\frac{s}{2} D^{2}}\partial_{\nu} (D^{2}) ||\, ||e^{(s-t) D^{2}}
Q_{\epsilon}||ds.  \end{array}$

\noindent
Now $||e^{(s-t) D^{2}} Q_{\epsilon}|| \leq 1$, and since
$\partial_{\nu} (D^{2})$ is a differential operator of order two with
bounded coefficients, $e^{-\frac{1}{2} D^{2}}\partial_{\nu} (D^{2})$
is a smoothing operator so has bounded norm.  Thus
$$
||e^{-\frac{s}{2} D^{2}}\partial_{\nu} (D^{2})|| \leq
||e^{-\frac{s-1}{2} D^{2}}||\,
||e^{-\frac{1}{2} D^{2}}\partial_{\nu} (D^{2})|| \leq
||e^{-\frac{1}{2}
D^{2}}\partial_{\nu} (D^{2})||
$$
for $t > 2$, as then $||e^{-\frac{s-1}{2} D^{2}}|| \leq 1$ for all $s
\geq t/2$.  Finally, $||Q_{\epsilon}e^{-\frac{s}{2} D^{2}}|| \leq
e^{-s\epsilon/2}$, so the last integral is bounded by a multiple of
$$
\epsilon^{-1}(e^{-(t\epsilon/4)} - e^{-(t\epsilon/2)}) =
t^{-\delta}(e^{-(t^{1+\delta}/4)} - e^{-(t^{1+\delta}/2)}) <
t^{-\delta}e^{-(t^{1+\delta}/4)}.
$$
This in turn is bounded by a multiple of $e^{-(t^{1+\delta}/8)}$, for
$t$ sufficiently large.

The change of variables $s\to t-s$ transforms the integral $\dd
\int^{t/2}_0 Q_{\epsilon} e^{-s D^{2}} \partial_{\nu} (D^{2}) e^{(s-t)
D^{2}} Q_{\epsilon}ds$ to the integral $\dd \int^t_{t/2} Q_{\epsilon}
e^{(s-t) D^{2}} \partial_{\nu} (D^{2}) e^{-s D^{2}} Q_{\epsilon}ds$,
so this satisfies the same estimate.  Replacing $D^{2}$ by $D^{2}/k$
then gives the estimate of the lemma.
\end{proof}
\begin{lemma}\label{estsqr} As $t \to \infty$,
$||Q_{\epsilon}\pa_{\nu}(e^{-tD^{2}}\sqrt{t}D)Q_{\epsilon}||$ is
bounded by a multiple of $e^{-(t^{1+\delta}/32)}$.
\end{lemma}

\begin{proof} Observe that
$$Q_{\epsilon}\pa_{\nu}(e^{-tD^{2}}\sqrt{t}D)Q_{\epsilon}=
Q_{\epsilon}\pa_{\nu}(e^{-tD^{2}/2}\sqrt{t}De^{-tD^{2}/2})
Q_{\epsilon}=
$$
$$
Q_{\epsilon}\pa_{\nu}(e^{-tD^{2}/2})\sqrt{t}De^{-tD^{2}/2}
Q_{\epsilon}+
Q_{\epsilon}e^{-tD^{2}/2}\pa_{\nu}(\sqrt{t}D)e^{-tD^{2}/2}
Q_{\epsilon}+
Q_{\epsilon}e^{-tD^{2}/2}\sqrt{t}D\pa_{\nu}(e^{-tD^{2}/2})
Q_{\epsilon}=
$$
$$
Q_{\epsilon}\pa_{\nu}(e^{-tD^{2}/2})Q_{\epsilon}\sqrt{t}De^{-tD^{2}/2}
Q_{\epsilon}+
Q_{\epsilon}e^{-tD^{2}/2}\pa_{\nu}(\sqrt{t}D)e^{-tD^{2}/2}
Q_{\epsilon}+
Q_{\epsilon}e^{-tD^{2}/2}\sqrt{t}DQ_{\epsilon}\pa_{\nu}(e^{-tD^{2}/2})
Q_{\epsilon}.
$$

The operators $De^{-tD^{2}/2}$, $\pa_{\nu}(D)e^{-tD^{2}/2}$
and $e^{-tD^{2}/2}D$ are all smoothing operators with norms bounded
independently of $t$, for $t$ large.  The fact that
$||Q_{\epsilon}e^{-tD^{2}/2}||  \leq e^{-t\epsilon/2} = e^{-(t^{1+\delta}/2)}$ 
and the estimate in Lemma~\ref{estk} give that
$||Q_{\epsilon}\pa_{\nu}(e^{-tD^{2}}\sqrt{t}D)Q_{\epsilon}||$ is
bounded by a multiple of
$\sqrt{t}(e^{-(t^{1+\delta}/16)}+e^{-(t^{1+\delta}/2)})$ which is
bounded by a multiple of $e^{-(t^{1+\delta}/32)}$, for $t$ large.
\end{proof}

\begin{lemma} As $t \to \infty$,
$||Q_{\epsilon}\pa_{\nu}(e^{-tD^{2}}{\dd \frac{I-e^{-tD^2}}{t D^2}}
\sqrt{t}D)Q_{\epsilon}||$ is bounded by a multiple of
$e^{-(t^{1+\delta}/32)}$.
\end{lemma}

\begin{proof}
$$
||Q_{\epsilon}\pa_{\nu}(e^{-tD^{2}}{\frac{I-e^{-tD^2}}{t D^2}}
\sqrt{t}D)Q_{\epsilon}|| =
||Q_{\epsilon}\pa_{\nu}(e^{-tD^{2}}\sqrt{t}D{\frac{I-e^{-tD^2}}{t D^2}}
)Q_{\epsilon}|| \leq
$$
$$
||Q_{\epsilon}\pa_{\nu}(e^{-tD^{2}}\sqrt{t}D)Q_{\epsilon}{\frac{I-e^{-tD^2}}{t
D^2}}
Q_{\epsilon}||+
||Q_{\epsilon}e^{-tD^{2}}\sqrt{t}D\pa_{\nu}({\frac{I-e^{-tD^2}}{t
D^2}}
)Q_{\epsilon}||
$$
and $||\dd{\frac{I-e^{-tD^2}}{t D^2}} || \leq 1$, so by Lemma~\ref{estsqr} the
first term immediately above satisfies the lemma.  If $G$  is the
Green's operator for $D$, the second term may be
written as
$Q_{\epsilon}(G/\sqrt{t})Q_{\epsilon}e^{-tD^{2}}tD^{2}\pa_{\nu}
(\dd{\frac{I-e^{-tD^2}}{tD^2}} )Q_{\epsilon},$
and $||Q_{\epsilon}G/\sqrt{t}|| \leq (t\epsilon)^{-1/2} =
t^{-(1+\delta)/2}$, which is bounded for $t$ large since $1+\delta
>0$.  The operator $tD^{2}{\dd \frac{I-e^{-tD^2}}{t D^2}} = I-e^{-tD^2}$, so
$$
tD^{2}\pa_{\nu}({\frac{I-e^{-tD^2}}{tD^2}} )=
-\pa_{\nu}(tD^{2}){\frac{I-e^{-tD^2}}{tD^2}} - \pa_{\nu}(e^{-tD^2}),
$$
and
$$
Q_{\epsilon}e^{-tD^{2}}tD^{2}\pa_{\nu}({\frac{I-e^{-tD^2}}{t
D^2}} )Q_{\epsilon}=
-Q_{\epsilon}e^{-tD^{2}}\pa_{\nu}(tD^{2}){\frac{I-e^{-tD^2}}{tD^2}}Q_{\epsilon}
- Q_{\epsilon}e^{-tD^{2}}\pa_{\nu}(e^{-tD^2})Q_{\epsilon}.
$$
Now,
$$
Q_{\epsilon}e^{-tD^{2}}\pa_{\nu}(tD^{2}) =
Q_{\epsilon}te^{-tD^{2}/2}e^{-tD^{2}/2}\pa_{\nu}(D^{2}),
$$
and $e^{-tD^{2}/2}\pa_{\nu}(D^{2})$ is a smoothing operator with norm
bounded independently of $t$, for $t$ large.
As
$$||Q_{\epsilon}te^{-tD^{2}/2}|| \leq te^{-t\epsilon/2} =
te^{-(t^{1+ \delta}/2)} < e^{-(t^{1+ \delta}/4)}
$$
for $t$ large, the term 
$\dd Q_{\epsilon}e^{-tD^{2}}\pa_{\nu}(tD^{2}){\frac{I-e^{-tD^2}}{tD^2}}Q_{\epsilon}$ 
 has norm bounded by a multiple of $e^{-(t^{1+\delta}/4)}$.
By Lemma~\ref{estk}, the  term 
$Q_{\epsilon}e^{-tD^{2}}\pa_{\nu}(e^{-tD^2})Q_{\epsilon} =
Q_{\epsilon}e^{-tD^{2}}Q_{\epsilon}\pa_{\nu}(e^{-tD^2})Q_{\epsilon}
$
is bounded by a multiple of
$e^{-(t^{1+\delta}/8)}$ (actually $e^{-t^{1+\delta}}$ if we use the
estimate $||Q_{\epsilon}e^{-tD^{2}}|| \leq e^{-t\epsilon} =
e^{-t^{1+\delta}}$).
\end{proof}

Thus we have the second inequality of Proposition \ref{estimates}.   
The third estimate follows immediately from the fact that both $P_t$ and $P_{\epsilon}$ are bounded.

\begin{lemma}\label{estep1} 
$||P_{\epsilon}\pa_{\nu}(e^{-tD^{2}})P_{\epsilon}||$ is bounded by a
multiple of $t^{1+(\delta/2)}$.
\end{lemma}

Note that $1+(\delta/2) > 1/2$, but by choosing $\delta$ close to
$-1$,
we can make $1+(\delta/2)$ as close to $1/2$ as we please.

\begin{proof}
$$
P_{\epsilon}\pa_{\nu}(e^{-tD^{2}})P_{\epsilon} =
-\int^t_0 P_{\epsilon} e^{-s D^{2}} \partial_{\nu} (D^{2}) e^{(s-t)
D^{2}} P_{\epsilon}ds =
-\int^t_0 P_{\epsilon} e^{-s D^{2}} P_{\epsilon}
[\partial_{\nu}(D)D + D\partial_{\nu}(D)]P_{\epsilon} e^{(s-t) D^{2}}
P_{\epsilon}ds.
$$
As $P_{\epsilon}$ is a smoothing operator, so are
 $P_{\epsilon}\pa_{\nu}(D)$ and $\pa_{\nu}(D)P_{\epsilon}$, since
$\pa_{\nu}(D)$ is a differential operator of order one with bounded
coefficients.  Since $\epsilon \to 0$ as $t \to \infty$, their norms are
bounded independently of $t$ for $t$ large.  Both $||P_{\epsilon}
e^{-s D^{2}}||$ and $|| e^{(s-t) D^{2}}P_{\epsilon}||$ are bounded by
$1$, and both $||P_{\epsilon} D||$ and $||D P_{\epsilon}||$ are bounded
by $\sqrt{\epsilon}$.  Thus
$||P_{\epsilon}\pa_{\nu}(e^{-tD^{2}})P_{\epsilon}||$ is bounded by a
multiple of $\dd\int^t_0 \sqrt{\epsilon}ds = \sqrt{\epsilon}t =
t^{1+(\delta/2)}.$
\end{proof}

\begin{lemma}\label{estep2}
$||P_{\epsilon}\pa_{\nu}(e^{-tD^{2}}\sqrt{t}D)P_{\epsilon}||$ is
bounded by a multiple of $t^{(3/2)+\delta}$.
\end{lemma}
Again note that we can make $3/2+\delta$ as close to $1/2$ as we
please.
\begin{proof}\

$\begin{array}{lll}||P_{\epsilon}\pa_{\nu}(e^{-tD^{2}}\sqrt{t}D)P_{\epsilon}||
& \leq &
||P_{\epsilon}\pa_{\nu}(e^{-tD^{2}})P_{\epsilon}\sqrt{t}DP_{\epsilon}||
+
||P_{\epsilon}e^{-tD^{2}}P_{\epsilon}\pa_{\nu}(\sqrt{t}D)P_{\epsilon}||\\
& \leq &
||P_{\epsilon}\pa_{\nu}(e^{-tD^{2}})P_{\epsilon}||\,\sqrt{t}\,
||DP_{\epsilon}||
+ \sqrt{t}\, ||P_{\epsilon}e^{-tD^{2}}P_{\epsilon}|| \, ||\pa_{\nu}(D)
P_{\epsilon}|| \\
& \leq & C_{1}(t^{1+(\delta/2)})\sqrt{t\epsilon} + C_{2}t^{1/2} =
C_{1}t^{(3/2)+\delta} + C_{2}t^{1/2} \leq Ct^{(3/2)+\delta}.  \end{array}$

\end{proof}

\begin{lemma}\label{estep3}
$||P_{\epsilon}\pa_{\nu}(e^{-tD^{2}}\dd{\frac{I-e^{-tD^2}}{t D^2}}
\sqrt{t}D)P_{\epsilon}||$ is bounded by a multiple of
$t^{(3/2)+\delta}$.
\end{lemma}

\begin{proof}
As $P_{\epsilon}\pa_{\nu}(e^{-tD^{2}}\dd{\frac{I-e^{-tD^2}}{t D^2}}
\sqrt{t}D)P_{\epsilon} = 
P_{\epsilon}\pa_{\nu}(e^{-tD^{2}}\sqrt{t}D\dd{\frac{I-e^{-tD^2}}{t D^2}}
)P_{\epsilon}$,   we have that
$$
P_{\epsilon}\pa_{\nu}(e^{-tD^{2}}\dd{\frac{I-e^{-tD^2}}{t D^2}}
\sqrt{t}D)P_{\epsilon} = 
P_{\epsilon}\pa_{\nu}(e^{-tD^{2}}\sqrt{t}D)P_{\epsilon}\dd{\frac{I-e^{-tD^2}}{t D^2}}
P_{\epsilon}  +
P_{\epsilon} e^{-tD^{2}}\sqrt{t}DP_{\epsilon}\dd\pa_{\nu}({\frac{I-e^{-tD^2}}{t D^2}}
)P_{\epsilon}.
$$
By the Lemma \ref{estep2} and the fact that
$||\dd {\frac{I-e^{-tD^2}}{t D^2}}|| \leq 1$, we need only consider the term
$$
||P_{\epsilon}e^{-tD^{2}}\sqrt{t}DP_{\epsilon}\pa_{\nu}({\frac{I-e^{-tD^2}}{t
D^2}})P_{\epsilon}||\leq
||e^{-tD^{2}}\sqrt{t}D P_{\epsilon} || \,
||P_{\epsilon}\pa_{\nu}({\frac{I-e^{-tD^2}}{t D^2}})P_{\epsilon}||
\leq
$$
$$
C\sqrt{t\epsilon} \,  ||P_{\epsilon}\pa_{\nu}({\frac{I-e^{-tD^2}}{t D^2}}
)P_{\epsilon}|| = Ct^{(1+\delta)/2}
||P_{\epsilon}\pa_{\nu}({\frac{I-e^{-tD^2}}{t D^2}})P_{\epsilon}||.
$$
Thus we need only show that
$||P_{\epsilon}\pa_{\nu}(\dd {\frac{I-e^{-tD^2}}{t
D^2}})P_{\epsilon}||$ is bounded by a multiple of $t^{1+(\delta/2)}$.
Note that
$$
\frac{d}{dr} (\frac{I-e^{-rD^2}}{ D^2}) = e^{-rD^2}
$$
so
$$\frac{d}{dr}(\pa_{\nu}(\frac{I-e^{-rD^2}}{ D^2}))=
\pa_{\nu}(\frac{d}{dr}(\frac{I-e^{-rD^2}}{ D^2})) =
\pa_{\nu}(e^{-rD^2}) = -\int^r_0 e^{-s D^{2}} \partial_{\nu} (D^{2})
e^{(s-r) D^{2}} ds.
$$
Thus
$$
\pa_{\nu}(\frac{I-e^{-tD^2}}{ D^2})= \int_{0}^{t}
\frac{d}{dr} (\pa_{\nu}(\frac{I-e^{-rD^2}}{ D^2}))dr  =
-\int_{0}^{t} \int_{0}^{r} e^{-s D^{2}} \partial_{\nu}
(D^{2}) e^{(s-r) D^{2}} ds\, dr,
$$
and
\begin{multline*}
||P_{\epsilon}\pa_{\nu}({\frac{I-e^{-tD^2}}{t D^2}})P_{\epsilon}|| =
||\frac{1}{t} P_{\epsilon}\pa_{\nu}(\frac{I-e^{-tD^2}}{
D^2})P_{\epsilon}||=
||\frac{1}{t}\int_{0}^{t} \int_{0}^{r} e^{-s D^{2}}  P_{\epsilon}
\partial_{\nu} (D^{2})P_{\epsilon} e^{(s-r) D^{2}} ds\,
dr|| \leq
\\ 
\frac{1}{t}\int_{0}^{t} \int_{0}^{r} ||e^{-s D^{2}}|| \, ||
P_{\epsilon}[
\partial_{\nu}(D)D + D\partial_{\nu}(D)]P_{\epsilon}|| \, || e^{(s-r)
D^{2}}|| ds\, dr \leq 
\frac{1}{t}\int_{0}^{t} \int_{0}^{r} C \sqrt{\epsilon} \, ds\, dr =
Ct^{1+(\delta/2)}.
\end{multline*}
\end{proof}
This finishes the proof of Proposition \ref{estimates}
\end{proof}

To finish the proof of Theorem \ref{limit}, first note that the estimates of Proposition \ref{estimates} remain true with $\pa_{\nu}$ replaced by $\delta$.  This follows from the fact that  for $T \in A_\mfS (\cG,\cS \otimes  E)\subset \wcA_{\mfS}$, $\delta T$ involves only $T$ and $\pa_{\nu} T$.  Similarly, $\delta P_0$,  $\delta (\alpha P_0)$, $\delta P_{\epsilon}$,  and $\delta Q_{\epsilon}$ are bounded operators.    Finally, for $\wT_1, \wT_2 \in  \wcA_{\mfS}$,  $\wT_1 * \wT_2 = \wT_1 \Theta \wT_2$.  But multiplication by 
$ \Theta=\left( \begin{array}{cc}1 &\quad 0\\  0 & \quad \theta \end{array} \right) $ 
is a bounded operation, so we may ignore it  with impunity in the norm estimates of products below.

Since   $P_{t} = \alpha P_{0} + P_{\epsilon}P_{t}P_{\epsilon} +
Q_{\epsilon}P_{t}Q_{\epsilon}$,
$$
\tr(P_{t}(\delta P_{t})^{2k}) = \tr (\alpha P_{0}(\delta (\alpha P_{0}))^{2k}) +
\tr(\alpha P_{0}(\delta P_{t})^{2k} - \alpha P_{0}(\delta (\alpha P_{0}))^{2k}) +
\tr(P_{\epsilon}P_{t}P_{\epsilon}(\delta P_{t})^{2k}) +
\tr(Q_{\epsilon}P_{t}Q_{\epsilon}(\delta P_{t})^{2k}).
$$
For any integer $\ell \geq 0$,
$$
||D^{2\ell}Q_{\epsilon}P_{t}Q_{\epsilon}(\delta  P_{t})^{2k}|| =
||D^{2\ell}Q_{\epsilon}e^{-tD^{2}/4} Q_{\epsilon}
\what{P}_{t}(\delta  P_{t})^{2k}|| \leq
||D^{2\ell}Q_{\epsilon}e^{-tD^{2}/4} Q_{\epsilon}|| \,
||\what{P}_{t}(\delta  P_{t})^{2k}||.
$$
Now
$$
\delta (P_{t}) = \delta (\alpha P_{0}) +
\delta (P_{\epsilon}P_{t}P_{\epsilon}) +
\delta (Q_{\epsilon}P_{t}Q_{\epsilon}),
$$
and $||\what{P}_{t}||$ is bounded independently of $t$.  So  $||\what{P}_{t}(\delta P_{t})^{2k}||$ is bounded by a multiple of
$$
||(\delta  P_{t})^{2k}|| =
||\delta (\alpha P_{0}) +
\delta (P_{\epsilon}P_{t}P_{\epsilon}) +
\delta (Q_{\epsilon}P_{t}Q_{\epsilon})||^{2k}
\leq Ct^{2k(\frac{1}{2} + a )}
$$
where $a > 0$ is a number to be chosen later (as close to zero as we
please).
On the other hand, for $t$ sufficiently large (so that $t^{1+\delta} > 4\ell$), the maximum of
$z^{\ell}e^{-tz/4}$ on the interval $[\epsilon, \infty)$ occurs at
$\epsilon$, so
$$
||D^{2\ell}Q_{\epsilon}e^{-tD^{2}/4}Q_{\epsilon}|| \leq
\epsilon^{\ell}e^{-t\epsilon/4} =
t^{\delta\ell}e^{-(t^{(1+\delta)}/4)}
$$
so
$$
||D^{2\ell}Q_{\epsilon}P_{t}Q_{\epsilon}(\delta P_{t})^{2k}|| \leq
Ct^{2k(\frac{1}{2}+a)}t^{\delta \ell}e^{-(t^{(1+\delta)}/4)}
$$
which goes to zero as $t \to \infty$.
The proof of Theorem 2.3.13 of \cite{H-L:1990} shows that this
implies that
$\tr(Q_{\epsilon}P_{t}Q_{\epsilon}(\delta P_{t})^{2k})$ is
pointwise bounded on $M$ and converges pointwise to zero as $t \to
\infty$.  As $\Phi$ is integration over a compact set, the bounded
convergence theorem gives
\begin{Equation}\label{eqinf1} \hspace{2in}
$ \dd \lim_{t \to \infty}\Phi\circ \tr(Q_{\epsilon}P_{t}
Q_{\epsilon}(\delta P_{t})^{2k}) = 0.$
\end{Equation}

Now consider
$\Phi \circ \tr(P_{\epsilon}P_{t}P_{\epsilon}(\delta P_{t})^{2k}) =
\Phi \circ\tr(P_{\epsilon}(\delta P_{t})^{2k}P_{\epsilon}
P_{t}P_{\epsilon})$.
The proof of Proposition 12 of \cite{H-L:1999} shows that
$$
||\Phi \circ
\tr(P_{\epsilon}(\delta P_{t})^{2k}P_{\epsilon}P_{t}P_{\epsilon})||_{T}
\,\, \leq \,\,
C||P_{\epsilon}(\delta P_{t})^{2k}P_{\epsilon}P_{t}|| \,
||\Phi \circ\tr(P_{\epsilon})||_{T}.
$$
Since $\Phi \circ\tr(P_{\epsilon})$ is $\cO(\epsilon^{\beta})$,
this is bounded by a multiple of
$$
t^{2k(\frac{1}{2}+a)}\epsilon^{\beta}\,\,  =\,\,
t^{2k(\frac{1}{2}+a)}t^{\delta\beta}\,\, = \,\,
t^{2k(\frac{1}{2}+a)+\delta\beta}.
$$
Recall that $-1 < \delta < -k/\beta$ and therefore we can choose $a > 0 $ so small
that
$$
2k(\frac{1}{2}+a)+\delta\beta < 0.
$$
Then
\begin{Equation}\label{eqinf2} \hspace{2in}
$ \dd \lim_{t \to \infty}\Phi \circ \tr(P_{\epsilon}P_{t}
P_{\epsilon}(\pa_{\nu}P_{t})^{2k}) = 0.$
\end{Equation}

Finally, consider the individual terms of
$\tr(\alpha P_{0}(\delta  P_{t})^{2k} - \alpha P_{0}(\delta  P_{0})^{2k})=
\tr(P_{0} \alpha (\delta  P_{t})^{2k} - P_{0}  \alpha (\delta  P_{0})^{2k}).$
Suppose the term $P_{0}A$ contains a
$\delta  (Q_{\epsilon}P_{t}Q_{\epsilon})$.  Then
$$
||P_{0}A|| \leq ||P_{0}|| \, ||\alpha||\,||\delta   (Q_{\epsilon}P_{t}Q_{\epsilon})
||^{\mu}||\delta  (\alpha P_{0})||^{\beta}||\delta  
(P_{\epsilon}P_{t}P_{\epsilon})||^{\gamma}
$$
where $\mu+ \beta + \gamma = 2k$ and $\mu > 0$.  Since $||\alpha||$ is bounded,  Proposition \ref{estimates}, gives that as $t \to \infty$,
$$
||P_{0}A|| \leq Ce^{-\mu(t^{1+\delta}/32)}t^{\gamma(\frac{1}{2}+a)}.
$$
For every positive integer $\ell$, $D^{2\ell}P_{0} = 0$,
so for every integer $\ell \geq 0$,
$||D^{2\ell}P_{0}A|| \to 0$ as $t \to \infty$.
Proceeding as in the proof of Equation \ref{eqinf1}, we have
$$
\lim_{t \to \infty}\Phi \circ \tr(P_{0}A) = 0.
$$

Now suppose that we have one of the remaining terms.  It must contain
a term of the form $\delta (P_{\epsilon}P_{t}P_{\epsilon})$.  As
$P_{\epsilon}^{2} = P_{\epsilon}$ and $\delta$ is a derivation, we may replace
$\delta  (P_{\epsilon}P_{t}P_{\epsilon})$ by
$$
\delta  (P_{\epsilon}^{2}P_{t}P_{\epsilon}) =
\delta  (P_{\epsilon})(P_{\epsilon}P_{t}P_{\epsilon}) +
P_{\epsilon}\delta  (P_{\epsilon}P_{t}P_{\epsilon}) =
\delta  (P_{\epsilon})(P_{\epsilon}P_{t}P_{\epsilon})P_{\epsilon} +
P_{\epsilon}\delta  (P_{\epsilon}P_{t}P_{\epsilon}).
$$
Using the
trace property of $\Phi \circ \tr$, we get two terms of the form $\Phi
\circ \tr(AP_{\epsilon})$.    As above, the proof of Proposition 12 of \cite{H-L:1999} shows that
$$
||\Phi \circ\tr(AP_{\epsilon})||_{T}  \,\, \leq \,\,  C||A|| \, ||\Phi\circ\tr(P_{\epsilon})||_{T}.
$$
Now $A$ is a product of terms of the form $\alpha, P_{0}, \delta  (\alpha P_{0}), P_{\epsilon}, \delta  (P_{\epsilon}), P_{\epsilon}P_{t}P_{\epsilon},$ and   $\delta  (P_{\epsilon}P_{t}P_{\epsilon})$.  Each of these is bounded in norm, except the last which has norm bounded by a multiple of $t^{(\frac{1}{2}+a)}$.  As $A$ can contain no more that $2k$ terms of the form $\delta  (P_{\epsilon}P_{t}P_{\epsilon})$, and $\Phi \circ\tr(P_{\epsilon})$ is $\cO(\epsilon^{\beta})$, we have that  $||\Phi \circ\tr(AP_{\epsilon})||_{T}$ 
is bounded by a multiple of 
$$
t^{2k(\frac{1}{2}+a)}\epsilon^{\beta}\,\, =\,\,
t^{2k(\frac{1}{2}+a)}t^{\delta\beta}\,\, = \,\,
t^{2k(\frac{1}{2}+a)+\delta\beta}.
$$
By our choice of a, we have that the limit as $t \to \infty$ of these
terms is zero, just as in the proof of Equation \ref{eqinf2}.

This completes the proof of Theorem \ref{limit-k}

\section{Bismut superconnections}\label{bismut}

As noted above, in \cite{BH-I} we proved that the Chern character $\ch_a$ composed with the topological and analytic index maps of Connes-Skandalis \cite{CS:1984} yield the same map.  In particular,  for any  Dirac operator $D$, the Chern character of the topological index of $D$, coincides with the 
Chern character of the analytic index of $D$, i.e.
$$
\ch_{a}(\Ind_{t}(D^+)) = \ch_{a}(\Ind_{a}(D^+)).
$$

In \cite{H-L:1999}, we proved that $\ch_{a}(\Ind_{t}(D^+))$ is equal to
the Chern character of the index bundle of $D$ in another sense.  We
defined a ``connection'' $\nabla$ on the index bundle $[P_{0}]$ of $D$,
and defined the Chern character of $[P_{0}]$ to be the Haefliger class of
$\Tr(\alpha e^{-(\nabla ^2/2i\pi)})$.  We then used a Bismut
superconnection for foliations, \cite{Heitsch:1995}, to show that
$\ch_{a}(\Ind_{t}(D))$ contains the Haefliger form $\Tr (\alpha
e^{-(\nabla ^2/2i\pi)})$, provided that the assumptions of
Theorem~\ref{limit} are satisfied, but with the stronger assumption that the
Novikov-Shubin invariants of $D$ are greater than {\bf three} times
the codimension of $F$.  We will now show that whenever $P_{0}$ is
smooth, $\ch_{a}([P_{0}])$ contains the Haefliger form $\Tr ( \alpha
e^{-(\nabla ^2/2i\pi)})$, so the two definitions of the Chern
character of $[P_{0}]$ agree.

We first recall the construction of Bismut superconnections for $D$.  See
\cite{Bismut:1986}, \cite{BV:1987}, and also \cite{Heitsch:1995}.  Let
$\nabla^B$ be a Bott connection on $\nu^*_{s}$.  If $\omega_1,\ldots,
\omega_n$ is a local framing for $\nu^*_{s}$, then
$\nabla^B\omega_i=\sum^n_{j=1}\omega_j\otimes \theta^i_j$ where
$\theta^i_j$ are local one forms on $\cG$ and the $\theta^i_j$ satisfy
$d\omega_i=\dd \sum\limits^n_{i=1}\omega_j\land \theta^i_j$.  that is,
the composition $$ C^{\infty}(\nu^*_{s})
\stackrel{\nabla^B}{\rightarrow} C^{\infty}(\nu^*_{s}\otimes T^*\cG)
\stackrel{\land}{\rightarrow} C^{\infty}(\nu^*_{s} \land T^*\cG) $$
is just $\omega \to d\omega$. $\nabla^B$ induces a connection on
$\land\nu^*_{s}$ also denoted $\nabla^B$ so that
$$
C^{\infty}(\land\nu^*_{s}) \stackrel{\nabla^B}{\rightarrow}
C^{\infty}(\land\nu^*_{s} \otimes T^*\cG)
\stackrel{\land}{\rightarrow}
C^{\infty}(\land \nu^*_{s} \land T^*\cG) $$
is also just $\omega\to d\omega$.

Set $\cV=TF_{s}\oplus \nu_{s}\oplus\nu^*_{s}=T\cG\oplus\nu^*_{s}$ over
$\cG$, and define a symmetric bilinear form $g$ on $\cV$ as follows.
$TF_{s}$ and $\nu_{s}\oplus\nu^*_{s}$ are orthogonal and $g|TF_{s}$ is
$g_0|TF_{s}$.  The form $g|\nu_{s}\oplus\nu^*_{s}$ is given by the
canonical duality, i.e.  $\nu_{s}$ and $\nu^*_{s}$ are totally
isotropic and $g(X,\omega )=\omega (X)$ for $X\in\nu_{s}$, $\omega \in
\nu^*_{s}$.  In \cite{BV:1987}, p.~455.  it is shown that there is a
unique connection $\nabla$, the Bismut connection, on $\cV$ so that
$\nabla$ preserves $\nu^*_{s}$ and $g$, $\nabla|\nu^*_{s}=\nabla^B$
and for all $X$, $Y\in C^{\infty}(T\cG)$, $\nabla_XY-\nabla_YX=[X,Y]$.
Note that in general $\nabla$ does not preserve $T\cG$ but that for
$X,Y\in C^{\infty}(T\cG)$, $\nabla_XY-\nabla_Y X\in C^{\infty}(T\cG)$.

Consider the vector space $V=\R^p\oplus\R^n\oplus\R^{n*}$.  Define
a bilinear form $Q$ on $V$ as $g$ was on $\cV$, i.e.  $\R^p$ is
orthogonal to $\R^n\oplus \R^{n*}$, $Q|\R^p$ is the usual inner
product, and $Q|\R^n\oplus \R^{n*}$ is given by the canonical
duality.  Let $C(V,Q)$ be the associated Clifford algebra and set
$S_{0}=\land \R^{n*}\otimes S$ where $S$ is the spinor space
for $\R^p$ with the usual inner product.  Let $\rho$ be the
representation of the Clifford algebra of $\R^p$ in $S$.  Then
$S_{0}$ is the spinor space for $C(V,Q)$ with the Clifford
multiplication being defined by $$\begin{array}{rcl} \rho_0(X)(\omega
\otimes s) &= &(-1)^{\deg \omega}\omega\otimes \rho(X)s\\
\rho_0(Y)(\omega\otimes s)    &= &-2i(Y)\omega \otimes s\\
\rho_0(\phi)(\omega\otimes s) &=& \phi\land\omega\otimes s
\end{array}$$
\noindent
for $X\in \R^p$, $Y\in\R^n$, $\phi\in\R^{n*}$, $\omega\in \wedge
\R^{n*}$, $s \in S$.  See \cite{BV:1987}, p.~456 and for general
facts
about spinors and Clifford algebras, \cite{LawsonMichelson}.

The above fact allows Berline and Vergne to give a beautiful and
concise definition of Bismut superconnections for fiber
bundles which was extended to foliations in \cite{Heitsch:1995}.  Recall that $\cS$ is the spinor bundle along the leaves of $F_{s}$, and 
consider the vector bundle $\cS_{0}=\land\nu^*_{s}\otimes \cS$ over
$\cG$ and the bundle of Clifford algebras $C(\cV)$ over $\cG$
associated to $\cV,g$.  Then $\cS_{0,y}$, the fiber over $y \in \cG$
of $\cS_{0}$, is a module for the algebra $C(\cV)_y$ and we denote
the module action also by $\rho_0$.  The connection $\nabla$ on $\cV$
induces a connection $\nabla$ on $\cS_{0}$ (\cite{BV:1987}, p.~456;
or more generally \cite{LawsonMichelson}, Ch.~4).  Let $E$ be a vector
bundle with connection over $\cG$ as in Section~\ref{char}.  We shall also denote by $\nabla$
the tensor product connection on $\cS_{0}\otimes E$.

A Bismut superconnection $\mathbb B$ for $F_{s}$ and $E$ is the Dirac
type operator on
$ C^{\infty}_{c}(\cS_{0}\otimes E)$ defined as follows.  Let
$X_1,\ldots,
X_p$ be a local oriented orthonormal basis of $TF_{s}$, and
$X_{p+1},\ldots, X_{p+n}$ a local basis of $\nu_{s}$.  Let
$X^*_1,\ldots,
X^*_{p+n}$ be the dual basis in $TF_{s}\oplus \nu^*_{s}$, i.e.
$X^*_i=X_i$ for
$1\le i\le p$, $X^*_i=\omega_i$, for $\ p+1\le i\le p+n$ where
$\omega_i\in
\nu^*_{s}$ and $\omega_i(X_j)=\delta_{ij}$.  Set
$$
\mathbb B =\sum^{p+n}_{i=1} \bigl(\rho_0(X^*_i)\otimes
1\bigr)\nabla_{X_i}
=\sum^p_{i=1}\rho(X_i)\nabla_{X_i} +
\sum^{p+n}_{i=p+1}\omega_i\nabla_{X_i}.
$$
$\mathbb B$ does not depend on the choice of $X_1,\ldots, X_{p+n}$.

Since $\cS=\cS^+\oplus \cS^-$ is $\Z_2$ graded and
$\land\nu^*_{s}$ is $\Z$ graded, $\cS_0=\land\nu^*_{s}\otimes
\cS$ has a total $\Z_2$ grading and we write $\cS_0 = \cS_0^+\oplus
\cS_0^-$.  We then have an associated $\Z_2$ grading $\cS_0 \otimes E
= (\cS_0^+ \otimes E) \oplus (\cS_0^- \otimes E)$.  It is immediate
from the fact that $\nabla$ preserves the grading that $\mathbb B$ is
an odd operator, i.e.  $\mathbb B$ maps $C^{\infty}(\cS_0^+\otimes E)$
to $C^{\infty}(\cS_0^-\otimes E)$ and vice-versa.

Finally, we may use the $\Z$ grading on $\land\nu^*_{s}$ to grade the
operator $\mathbb B$, i.e.  $\mathbb B=\mathbb B^{\, [0]} +\mathbb
B^{\, [1]}+\cdots$ where $\mathbb B^{\, [i]}:C^{\infty} (\land^k\nu^*_{s}\otimes
\cS \otimes E) \to C^{\infty} (\land^{k+i}\nu^*_{s}\otimes \cS \otimes
E)$.

It is straightforward to check that
\begin{proposition}
The term $\mathbb B^{\, [1]}$ is a quasi-connection $\nabla^\nu$ for $E
\otimes \land \nu^{*}_{s}$as defined in Section \ref{chern}.
\end{proposition}

Recall, \cite{H-L:1999}, that a {\em connection} on the {\em index
bundle} of $D$ is defined
by
$$
\nabla = P_{0} \mathbb B^{\, [1]} P_{0}.
$$
 For this to be well defined, we must require that $P_{0}$
is smooth.

\begin{theorem}
Suppose that $P_{0}$ is smooth.  Then  $\ch_{a}([P_{0}])$ contains the
Haefliger form
$\Tr ( \alpha e^{-(\nabla ^2/2i\pi)})$
\end{theorem}

\begin{proof}\
First we calculate $\nabla ^2$.
\begin{eqnarray*}
\nabla ^2 & = & P_{0} \mathbb B^{\, [1]} P_{0} \mathbb B^{\, [1]} P_{0} \\
& = & P_{0} [\mathbb B^{\, [1]},  P_{0}] \mathbb B^{\, [1]}  P_{0}  + P_{0}
(\mathbb
B^{\, [1]})^2 P_{0} \\
& = & P_{0} [\mathbb B^{\, [1]}, P_{0}] [\mathbb B^{\, [1]}, P_{0}] + P_{0}
[\mathbb B^{\, [1]}, P_{0}] P_{0} \mathbb B^{\, [1]} + P_{0} (\mathbb
B^{\, [1]})^2 P_{0} \\
& = & P_{0} [\mathbb B^{\, [1]}, P_{0}] [\mathbb B^{\, [1]}, P_{0}] + P_{0}
(\mathbb B^{\, [1]})^2 P_{0}.
\end{eqnarray*}
The last equality is a consequence of the relation $P_{0} [\mathbb
B^{\, [1]}, P_{0}] P_{0} = 0$ which is true since $P_{0}^2=P_{0}$ and since
$[\mathbb B^{\, [1]},\cdot]$ is a derivation.  This derivation is precisely
$\pa_\nu$, so $(\mathbb B^{\, [1]})^2 = \theta$ as in Section~\ref{chern}.  Thus
$$
\nabla ^2  =  P_{0} (\pa_\nu P_{0})^2 +  P_{0} \theta P_{0},
$$

and
$$
\nabla ^{2k} = (P_{0} (\pa_\nu P_{0})^2 +  P_{0} \theta P_{0})^{k}.
$$
Note that
$$
\pa_{\nu}(P_{0}) = \pa_{\nu}(P_{0}P_{0}) = \pa_{\nu}(P_{0})P_{0} +
P_{0}\pa_{\nu}(P_{0}),
$$
so
$$
\pa_{\nu}(P_{0})P_{0} = \pa_{\nu}(P_{0}) - P_{0}\pa_{\nu}(P_{0}).
$$
Using this twice, one can easily show that
$$
P_{0}\pa_{\nu}(P_{0})\pa_{\nu}(P_{0})=
P_{0}\pa_{\nu}(P_{0})\pa_{\nu}(P_{0})P_{0} .
$$
Then a simple induction argument shows that
$$
(P_{0} (\pa_\nu P_{0})^2 +  P_{0} \theta P_{0})^{k} = 
P_{0} ((\pa_\nu P_{0})^2 +  P_{0} \theta P_{0})^{k}
$$
Thus,
$$
\Tr ( \alpha \nabla^{2k}) = \Tr (\alpha P_{0} ((\pa_\nu P_{0})^2 +  P_{0} \theta P_{0})^{k}),
$$
and comparing with 
Equation \ref{PO}, we see that $\ch_{a}([P_{0}])$ contains the Haefliger form
$\Tr ( \alpha e^{-(\nabla ^2/2i\pi)})$.
\end{proof}


\begin{thebibliography}{BGV92}

\bibitem[A75]{Atiyah:1975} M.  F.  Atiyah, {\em Elliptic operators,
discrete
groups and von Neumann algebras},  Asterisque {\bf 32/33} (1976)
43-72.

\bibitem[BH-I]{BH-I} M-T. Benameur and J. Heitsch,  {\em
Index theory and Non-Commutative Geometry I. Higher Families Index Theory},
to appear in  {\em K-Theory} .

\bibitem[BH-III]{BH-III} M-T. Benameur and J. Heitsch,  {\em
Index theory and Non-Commutative Geometry III. The Haefliger  Higher Signatures for Foliations},  in preparation.

\bibitem[BV87]{BV:1987}  N. Berline and M. Vergne,  {\em  A proof of
Bismut local index theorem for a family of Dirac operators.} {\sl
Topology\/} {\bf 26} (1987), 435--463.

\bibitem[B86]{Bismut:1986} J.-M.  Bismut, {\em The Atiyah-Singer index
theorem for families of Dirac operators: two heat equation proofs. }
{\sl Invent.  Math.} {\bf 83} (1986), 91--151.

\bibitem[Con79]{Connes:1979} A. Connes, {\em Sur la th\'eorie de
l'int\'egration non commutative}.  Lect. Notes in Math. {\bf 725},
1979.

\bibitem[Con81]{Connes:1981} A.~Connes. {\em A survey of foliations
and operator algebras}.  Operator algebras and applications, Part I,
Proc.
Sympos.  Pure Math {\bf 38}, Amer.  Math.  Soc., (1982) 521-628.

\bibitem[Con85]{ConnesIHES} A.~Connes. {\em Noncommutative differential geometry.}  Inst. Hautes Etudes Sci. Publ. Math.  No. 62 (1985), 257--360.

\bibitem[Con86]{ConnesFundamental} A.~Connes. {\em Cyclic cohomology and the transverse fundamental class of a foliation.}  Geometric methods in operator algebras (Kyoto, 1983),  52-144, Pitman Res. Notes Math. Ser., 123, Longman Sci. Tech., Harlow, 1986. 

\bibitem[Con94]{ConnesBook} A.~Connes,
\newblock {\em Noncommutative Geometry}, Academic Press, New York,
1994.

\bibitem[CM91]{ConnesMoscovici} A.~Connes and H. Moscovici,
{\em Cyclic cohomology and the Novikov conjecture for hyperbolic
groups}, {\em Topology} {\bf 29} (1990) 345-388.

\bibitem[CS84]{CS:1984} A.~Connes, and G.~Skandalis.
{\em The longitudinal index theorem for foliations},
\newblock {\em Publ.  RIMS Kyoto} {\bf 20} (1984) 1139-1183.

\bibitem[CQ97]{CuntzQuillen} J. Cuntz and D. Quillen,
{\em  Excision in bivariant periodic cyclic cohomology.}  Invent. Math.  127  (1997),  no. 1.

\bibitem[GL03]{GL:2003} A.~Gorokhovsky, and J.~Lott.  {\em Local
index theory over \'{e}tale groupoids},
\newblock {\em J. Reine Angew. Math. 560, p. 151-198, 2003}.

\bibitem[GL05]{GL:2005} A.~Gorokhovsky, and J.~Lott.  {\em Local
index theory over \'{e}tale groupoids II},
\newblock {\em preprint}.

\bibitem[H80]{Hae:1980}A.~Haefliger.  {\em Some remarks on
foliations with minimal leaves},  J.  Diff.  Geo. {\bf
15} (1980) 269--284.

\bibitem[He95]{Heitsch:1995}J.~L. Heitsch.
{\em Bismut superconnections and the {C}hern character for
{D}irac operators on foliated manifolds},
\newblock {\em K-Theory} {\bf 9} (1995) 507--528.

\bibitem[HL90]{H-L:1990}J.~L. Heitsch and C.~Lazarov.
{\em A {L}efschetz theorem for foliated manifolds},
Topology, {\bf 29} (1990) 127--162.

\bibitem[HL99]{H-L:1999}
J.~L. Heitsch and C.~Lazarov.
{\em A general families index theorem},
K-Theory, {\bf 18} (1999) 181--202.

\bibitem[HL02]{H-L:2002}
J.~L. Heitsch and C.~Lazarov.
{\em Riemann-Roch-Grothendieck and torsion for foliations}
J.  Geo.  Anal. {\bf 12} (2002) 437--468.

\bibitem[LM89]{LawsonMichelson} H.  B. Lawson and  M.-L. Michelson,
{\em Spin geometry}, Princeton Math.  Series {\bf 38}, Princeton, 1989.

\bibitem[MN96]{MelroseNistor} R. Melrose and V. Nistor, {\em Homology of pseudodifferential operators I (manifolds with boundary)}, to appear in Amer. J. Math.

\bibitem[Nis93]{NistorInv} V. Nistor, {\em A bivariant Chern-Connes character.}  Ann. of Math. (2)  138  (1993),  no. 3, 555--590.

\bibitem[NWX96]{NistorWeinsteinXu} V. Nistor, A. Weinstein and P. Xu,
{\em Pseudodifferential operators on differential groupoids}  Pacific
J. Math.  \textbf{189} (1999),  no. 1, 117--152.

\bibitem[RS80]{Reed-Simon}
M.~Reed, and B.~Simon.  \newblock {\em Methods of Modern Mathematical
Physics I: Functional Analysis},  Academic Press, New York, 1980.

\end{thebibliography}
\end{document}